\documentclass[]{iopart}

\usepackage{t1enc}
\usepackage{graphicx}
\usepackage[english]{babel}
\usepackage[T1]{fontenc}
\usepackage[latin1]{inputenc}
\usepackage{url}
\usepackage{dsfont}
\usepackage{graphicx}
\usepackage{graphics}
\usepackage{subfig}
\usepackage{color}
\usepackage{iopams} 
\usepackage{setstack}
\usepackage{amstext}
\usepackage{amsopn}
\usepackage{amsthm}
\usepackage{enumitem}
\usepackage{siunitx}
\usepackage{hyperref}
\hypersetup{colorlinks=true, breaklinks=true, linkcolor=darkblue, urlcolor=darkblue,citecolor=darkblue}

\definecolor{darkblue}{rgb}{0,0,0.6}
\definecolor{darkgreen}{rgb}{0,0.6,0}
\definecolor{darkred}{rgb}{0.6,0,0}
\definecolor{darkorange}{rgb}{1,0.4,0}

\DeclareMathOperator*{\mydef}{\mathrel{\mathop:}=}

\newcommand{\R}{\mathbb{R}}
\newcommand{\J}{\mathcal{J}}

\newcommand{\norm}[1]{\Vert #1 \Vert}
\newcommand{\sqnorm}[1]{\Vert #1 \Vert_2^2}

\newcommand{\argminsub}[1]{\underset{{ #1 }}{{\rm argmin}}}

\newcommand{\termabb}[2]{\emph{#1} (\emph{#2})}

\graphicspath{{figures/}}

\begin{document}

\title{Accelerated High-Resolution Photoacoustic Tomography via Compressed Sensing}

\author{Simon Arridge$^1$, Paul Beard$^2$, Marta Betcke$^1$, Ben Cox$^2$, Nam Huynh$^2$, Felix Lucka$^1$, Olumide Ogunlade$^2$  and Edward Zhang$^2$}

\address{$^1$ Department of Computer Science, University College London, WC1E 6BT London, UK}
\address{$^2$ Department of Medical Physics and Bioengineering, University College London, WC1E 6BT London, UK}

\ead{f.lucka@ucl.ac.uk}

\begin{abstract}
Current 3D photoacoustic tomography (PAT) systems offer either high image quality or high frame rates but are not able to deliver high spatial and temporal resolution simultaneously, which limits their ability to image dynamic processes in living tissue (4D PAT). A particular example is the planar Fabry-P\'{e}rot (FP) photoacoustic scanner, which yields high-resolution 3D images but takes several minutes to sequentially map the incident photoacoustic field on the 2D sensor plane, point-by-point. However, as the spatio-temporal complexity of many absorbing tissue structures is rather low, the data recorded in such a conventional, regularly sampled fashion is often highly redundant. 
We demonstrate that combining model-based, variational image reconstruction methods using spatial sparsity constraints with the development of novel PAT acquisition systems capable of sub-sampling the acoustic wave field can dramatically increase the acquisition speed while maintaining a good spatial resolution: First, we describe and model two general spatial sub-sampling schemes. Then, we discuss how to implement them using the FP interferometer and demonstrate the potential of these novel compressed sensing PAT devices through simulated data from a realistic numerical phantom and through measured data from a dynamic experimental phantom as well as from in-vivo experiments. Our results show that images with good spatial resolution and contrast can be obtained from highly sub-sampled PAT data if variational image reconstruction techniques that describe the tissues structures with suitable sparsity-constraints are used. In particular, we examine the use of total variation (TV) regularization enhanced by Bregman iterations. These novel reconstruction strategies offer new opportunities to dramatically increase the acquisition speed of photoacoustic scanners that employ point-by-point sequential scanning as well as reducing the channel count of parallelized schemes that use detector arrays.
\end{abstract}

\submitto{\PMB}
\maketitle


\section{Introduction} \label{sec:Intro}

\termabb{Photoacoustic Tomography}{PAT} is an emerging biomedical, "Imaging from Coupled Physics"-technique \cite{ArSc12} based on laser-generated ultrasound (US). It allows the rich contrast afforded by optical absorption to be imaged with the high spatial resolution of ultrasound. Furthermore, the wavelength dependence of the optical absorption can, in principle, be utilized to provide spectroscopic (chemical) information on the absorbing molecules (chromophores). While PAT's potential for an increasing variety of clinical applications is currently being explored 
\cite{ZaVaGa14,TaNt15}, it is already widely used in preclinical studies to examine small animal anatomy, physiology and pathology \cite{YaWa14,XiWa14,JaLaOgTrCoZhJoPiPhMaLyPePuBe15}. 
For further applications and references we refer to the reviews \cite{Wa09,Bea11,NiCh14}.\\
To obtain high quality three-dimensional (3D) photoacoustic (PA) images with a spatial resolution around one hundred \si{\micro \meter}, acoustic waves with a frequency content typically in the range of tens of \si{\mega \hertz} need to be sampled over \si{\centi \meter} scale apertures. Hence, satisfying the spatial Nyquist criterion necessitates sampling intervals on a scale of tens of \si{\micro \meter}, which requires scanning several thousand detection points. Using sequential scanning schemes, such as the \termabb{Fabry-P\'{e}rot based PA scanner}{FB scanner} or mechanically scanned piezoelectric receivers, this inevitably results in long acquisition times. In principle, this can be overcome by using an array of detectors. However, a fully sampled array comprising several thousand elements, each with their own signal conditioning electronics and radio frequency analog-to-digital (RF A-D) electronics, would be prohibitively expensive, and efficient multiplexing is challenging due to the low pulse repetition frequency of current excitation lasers. The slow acquisition speed currently limits the use of PAT for applications where movement of the target will cause image artefacts and prohibits the examination of dynamic anatomical and physiological events in high resolution in real time, which is the goal in 4D PAT. \\
A different approach to accelerate sequential PAT scanners relies on the key observation that, in many situations, the spatial complexity of many of the absorbing tissue structures is rather low, and therefore, data recorded in a conventional, regularly sampled, fashion is highly redundant. It may be possible, therefore, to speed up the data acquisition without a significant loss of image quality by exploiting this redundancy and measuring a subset of the data chosen in such a way as to maximize its non-redundancy. This concept, established as the field of \termabb{compressed sensing}{CS} \cite{CaRoTa06,Do06,FoRa13}, has been applied to several imaging modalities with success, most notably to \termabb{magnetic resonance tomography}{MRI} \cite{LuDoPa07,TrDiAtAr14,HaFaBeTaSeGlHo15} and \termabb{computed tomography}{CT} \cite{RiBeFlKa11,HmKaKoLaNiSi13,JoSi15}. 
As the inverse problem of PAT is conceptually similar to these, there has been an increased interest in applying CS to PAT in different ways and for different types of scanners:
In \cite{PrLe09,GuLiSoWa10,ZhWaZh12,MeWaLiSo12,MeWaYiLiSo12}, 2D reconstructions from regularly sub-sampled circular and linear transducer arrays were computed using matrix-based algorithms and first promising results were obtained. An asymmetrical 2D circular sensor arrangement designed based on a low-resolution pre-image was examined in \cite{CaYuDuXuWaCaLi15}. In \cite{SaKrBeBuHa15}, the CS-based recovery of suitably transformed integrating line detector data followed by a universal backprojection in 2D was proposed. Another approach, using patterned excitation light, was presented in \cite{LiZhYi09,SuFeShShMaLiWu11}, which is a possible approach for photoacoustic microscopy, but is not applicable in PAT where light undergoes strong scattering. In our work, we further explore the potential for acceleration of high resolution 3D PAT using randomized, incoherent encoding of the PA signals and a sparsity-constrained image reconstruction via variational regularization enhanced by Bregman iterations. Extensive studies with simulated data from realistic numerical phantoms, measured data from experimental phantoms and in-vivo recordings are carried out to demonstrate which conditions have to be fulfilled to obtain high sub-sampling rates. To cope with the immense computational challenges, GPU computing is combined with matrix-free, state-of-the-art optimization approaches.\\ 
The remainder of the paper is organized as follows: Section 2 provides background on theoretical and practical aspects of PAT. Section \ref{sec:CSPAT} describes two general approaches to accelerate sequential scanning schemes by spatial sub-sampling / compressive sensing and exemplifies their implementation with an FP scanner. In Section \ref{sec:ImgRec}, we then describe how images can be reconstructed from the sub-sampled/compressed data. Our methods are exemplified by simulation studies in Section \ref{sec:SimStud} and by reconstruction from experimental data in Section \ref{sec:ExpData}. In Section \ref{sec:DisOutCon}, we summarize and discuss the results of our work and point to future directions of research. Table \ref{tbl:Abb} lists all commonly occurring abbreviations for later look-up.

\section{Background} \label{sec:Background}

\subsection{Basics of Photoacoustic Tomography} \label{subsec:BasicsPAT}

\setlength{\fboxsep}{3pt}

The photoacoustic (PA) signal is generated by the coupling of optical and acoustic wave propagation processes through the \emph{photoacoustic effect}: Firstly, the tissue is illuminated by a laser pulse with a duration of a few nanoseconds. Inside the tissue, the photons will be scattered and absorbed, the latter predominantly in regions with a high concentration of chromophores, such as haemoglobin. The photoacoustic effect occurs when a sufficient part of the absorbed optical energy is converted to heat (\emph{thermalised}) sufficiently fast and not re-emitted: The induced, local pressure increase $p_0$ initiates acoustic waves travelling in the tissue on the microsecond timescale. These waves can be measured by ultrasonic transducers at the boundary of the domain. \\
With several assumptions on the tissue's properties (see \cite{WaAn11} for a detailed discussion), the acoustic part of the signal generation can be modelled by the following 
initial value problem for the wave equation:
\begin{equation}
\fl \qquad (\partial_{tt} - c_0^2 \Delta) p(r,t) = 0, \qquad p(r,t = 0) = p_0, \qquad \partial_t p(r,t = 0) = 0. \label{eq:PATfwd}
\end{equation}
The measurement data $f$ consists of samples of $p(r,t)$ on the boundary of the domain. See \cite{Be11,LuRa13} for recent reviews on measurement systems. The computation of a high resolution reconstruction of $p_0$, usually referred to as the \termabb{photoacoustic image}{PA image}, from $f$ is the subject of this paper. Obtaining a high quality PA image is of crucial importance for any subsequent analysis, e.g., for \termabb{quantitative photoacoustic tomography}{QPAT} \cite{CoLaArBe12}, wherein the optical part of the signal generation is inverted based on the PA image.

\subsection{Nyquist Sampling in Space and Time} \label{subsec:BackNyquist}

Before we discuss how to sub-sample the incident photoacoustic field $p(r,t)$ in Section \ref{sec:CSPAT}, it is important to understand that a complete sampling requires a certain relation between sampling in $r$ and $t$: Imagine we measure the PA signal on the boundary of a domain with homogenous sound speed $c$ which is band-limited to $\omega_t^*$. Firstly, the Nyquist criterion requires us to sample the PA signal with a temporal spacing of $\delta_t < 1/(2 \omega^*_{t})$. Secondly, the PA signal is caused by incident acoustic waves coming from various spatial directions. As an incident wave with a wave vector $\mathbf{k}$ leads to a PA signal with frequency 
\begin{equation}
\omega_t = c \sqrt{k_x^2 + k_y^2 + k_z^2} = c \| \mathbf{k} \|_2,
\end{equation}
the band-limit $\omega_t^*$ of the PA signal implies that the spatial waves are limited by $\| \mathbf{k}^* \|_2 = \omega_t^*/c$. To resolve all the spatial information, the Nyquist criterion would require to sample the domain boundary with a spatial spacing of $\delta_r < 1/(2 \| \mathbf{k}^* \|_2) = c/(2 \omega_t^*)$.


\section{Compressed Photoacoustic Sensing} \label{sec:CSPAT}

\setlength{\fboxsep}{0.01\textwidth}
\begin{figure}[tb]
\centering
\subfloat[][\label{subfig:FPSL}]{\fbox{\includegraphics[width=0.475\textwidth]{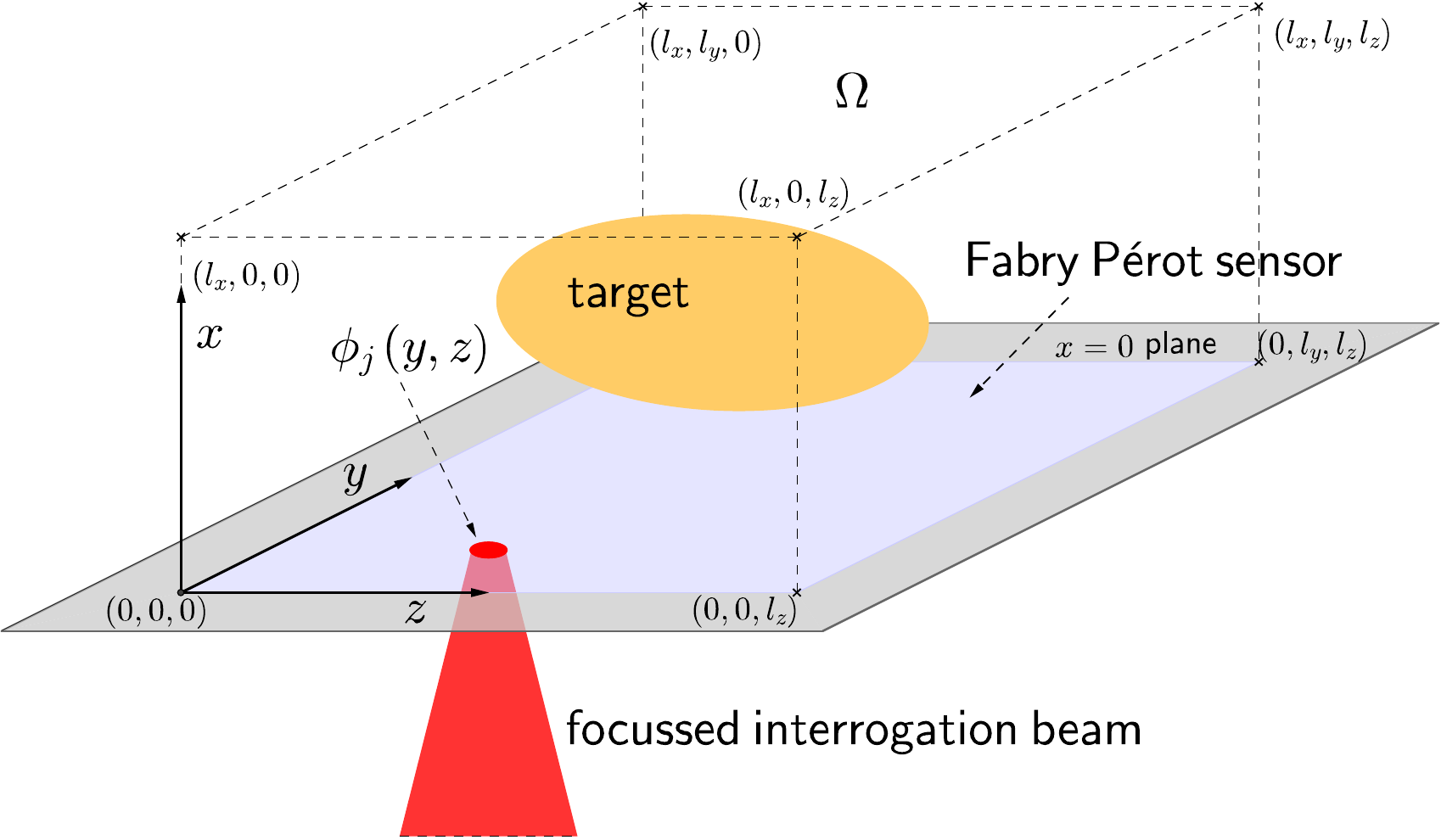}}}
\hfill
\subfloat[][\label{subfig:FPBP}]{\fbox{\includegraphics[width=0.475\textwidth]{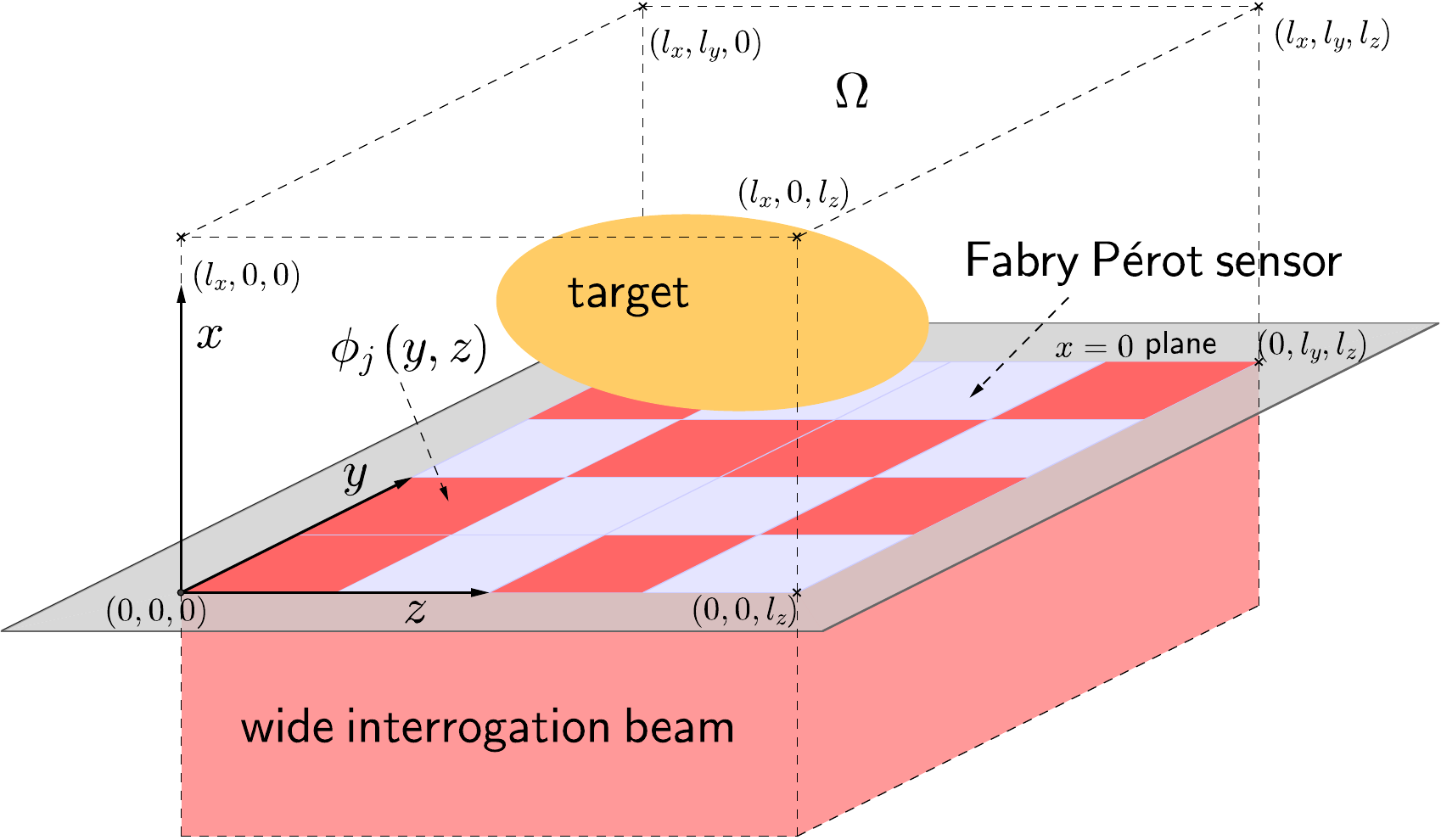}}}
\caption{Conceptual sketch of a Fabry-P\'{e}rot interferometer-based PAT setting with different spatial window functions $\phi_j(y,z)$: \protect\subref{subfig:FPSL} The interrogation beam is focused on a single location, leading to a very localized $\phi_j(y,z)$. \protect\subref{subfig:FPBP} A wide interrogation beam is used that has been pattered to produce a distributed, binary $\phi_j(y,z)$.}
   \label{fig:SubSam}
\end{figure}

\subsection{Sub-Sampling of the Photoacoustic Field}  \label{subsec:SubSam}

As discussed in Section \ref{sec:Intro}, a drawback of all current 3D PAT systems is that they offer either exquisite image quality or high frame rates but not both, partly due to the difficulty of realizing a scheme that would complete a scan with a sufficiently small $(\delta_r,\delta_t)$ in an acceptable acquisition time. In this section, we describe two novel sensing paradigms that aim to accelerate the data acquisition by spatially sub-sampling the incident photoacoustic field $p(r,t)$. In this study, the practical realization of these two approaches is achieved via the  \termabb{Fabry-P\'{e}rot}{FP} photoacoustic scanner as  described in Section \ref{subsec:CSFP}, but they are equally applicable to other sequential scanning systems, such as mechanically scanned piezoelectric detectors. \\
For simplicity but without loss of generality, we assume that a planar detection surface located at $x = 0$ is used, that a rectangular area $[0,l_y] \times [0,l_z]$ on it can be interrogated, and that the target is located in the region $\Omega = [0,l_x] \times [0,l_y] \times [0,l_z]$. The extension to other detection geometries is straight forward, but complicates the notation. The concrete measurement process can be modeled in the following way: The incident photoacoustic field on the the detection plane, $p(x=0,y,z,t)$, caused by the $j^{th}$ pulse of the excitation laser, is first multiplied by a spatial window function $\phi_j(y,z)$ and then integrated over the whole detection surface: 
\begin{equation}
f_j(t) = \int p(x=0,y,z, t) \phi_j(y,z) \; \rmd y \rmd z
 \label{eq:GenScan}
\end{equation}
The spatial sampling is followed by a temporal sampling, e.g., by measuring $f_j(t)$ at $t_i = i_t \delta_t$, $i_t = 1,\ldots,M_t$. The setting is sketched in Figure \ref{fig:SubSam}.

\subsubsection{Single-Point Scanning} \label{subsubsec:SiPo}

A standard approach is to try to focus $\phi_j(y,z)$ on a single point $(y_j,z_j)$ (cf. Figure \ref{subfig:FPSL}), which belongs to a set of $M = M_y \times M_z$ points forming a regular grid with spacings $\delta_y = l_y/M_y$ and $\delta_z = l_z/M_z$. We will refer to the data obtained by this scanning pattern as the \termabb{conventional}{cnv.} data. Ideally, these spacings are chosen fine enough to assure that the Nyquist criterion is satisfied in space and time (cf. Section \ref{subsec:BackNyquist}).\\
Due to the spatio-temporal characteristics of the wave propagation, the pressure time series recorded at two neighbouring locations on the regular grid provide very similar information. To reduce the coherence of the measured time series, we can instead also sample the photoacoustic field at a smaller number of randomly chosen points $(y_j,z_j), j=1,\ldots,M_c$, $M_c < M$, which would yield an acceleration factor of $M_{sub} = M/M_c$.  We will denote this sub-sampling strategy by \mbox{rSP-$M_{sub}$}. Choosing random points is firstly motivated by several results in compressed sensing theory \cite{FoRa13}, that point out the importance of randomness for designing sub-sampling pattern. Secondly, we want to avoid choosing sampling pattern that might lead to systematic biases. We will return to this point in the computational studies in Sections \ref{sec:SimStud} and \ref{sec:ExpData}.

\subsubsection{Patterned Interrogation Scanning}  \label{subsubsec:PatInter}

The second idea to accelerate the acquisition is to choose a series of orthogonal pattern ${\phi_j(y,z)}_j$ with each $\phi_j(y,z)$ being supported all over $[0,l_y] \times [0,l_z]$ with no particular focus on a single location (cf. Figure \ref{subfig:FPBP}). Again, choosing $M_c < M$ pattern would yield an acceleration factor of $M_{sub} = M/M_c$ over conventional single point scanning. As we will use a specific type of binary pattern based on scrambled Hadamard matrices described later (Section \ref{subsec:DiscFwdModel}), we will denote the sub-sampling strategy by \mbox{sHd-$M_{sub}$}.

\subsection{Implementation of Sub-Sampling Schemes using the FB Scanner}  \label{subsec:CSFP}

The planar, interferometry-based Fabry-P\'{e}rot photoacoustic scanner provides a convenient implementation of the sub-sampling strategies described above. In the standard set-up, the FP scanner performs the conventional single-point scanning described in Section \ref{subsubsec:SiPo}: For each pulse sent in by the excitation laser, the pressure time series at a different location on a grid is measured by interrogating the FP sensor head with an interrogation laser \cite{ZhLaBe08}. Due to the planar geometry and the option to introduce the excitation laser through the transparent detection plane, PA signals from a large range of anatomical targets can be scanned in the frequency range from DC to several tens of \si{\mega\hertz} on a scale of tens of \si{\micro\meter} ($\phi_j(y,z)$ in our model corresponds to the beam profile of the interrogation laser, cf. Figure \ref{subfig:FPSL}). Superficial features located a few \si{\milli \meter} below the skin surface can be imaged with high spatial resolution from a conventional FP scan (see, e.g., \cite{JaLaOgTrCoZhJoPiPhMaLyPePuBe15}). However, as described in Section \ref{sec:Intro},  obtaining such an image requires scanning several tens of thousands of locations which, due to the repetition rates of currently available excitation lasers, takes several minutes. \\
While the single point sub-sampling described in Section \ref{subsubsec:SiPo} is straight forward to implement with a standard FP scanner, implementing the patterned interrogation scheme requires several modifications: Instead of focusing the interrogation beam on a single location, the whole detection plane is illuminated  and the reflected beam is patterned before being focused into the photo diode. The spatial modulation is inspired by the working-principle of the celebrated "single-pixel Rice camera" \cite{DuDaTaLaTiKeBa08}: A \termabb{digital micromirror device}{DMD} is used to block rectangular sections of the reflected interrogation beam which creates binary pattern. Note that this hardware realization slightly differs from the conceptual sketch in Figure \ref{subfig:FPBP}, where the direct interrogation beam is patterned. Further details of the \textit{patterned-illumination} scanner can be found in \cite{HuZhBeArBeCo14,HuZhBeArBeCo15,HuZhBeArBeCo16}. 


\section{Image Reconstruction from Sub-Sampled Data} \label{sec:ImgRec}

In this section, we describe a model of the accelerated data acquisition and how to invert it.  

\subsection{Continuous Forward Model}  \label{subsec:FwdModel}

The PAT forward operator, $A$, maps a given initial pressure $p_0$ in the volume $[0,l_x] \times [0,l_y] \times [0,l_z]$ to the time dependent pressure on the detection plane, $\bar{p} \mydef p(x=0,y,z,t)$, as determined by \eref{eq:PATfwd}: $\bar{p} = A p_0$. A more detailed discussion of the operator $A$ and its adjoint can be found in \cite{ArBeCoLuTr16}. In this work, we assume that the target's dynamics are slow enough to be well approximated by a constant within the acquisition time. A sensing operator $C$ implements \eref{eq:GenScan} to produce the measured data $f^c \in \R^{M_c M_t}$ from the detection plane pressure $\bar{p}$:
\begin{equation}
f^c = C \bar{p} + \varepsilon = C A p_{0} + \varepsilon, \label{eq:Fwd}
\end{equation}
where we added the term $\varepsilon$ to account for additive measurement noise. In the following, $C$ denotes the sampling operator: The conventional point-by-point scan described in Section \ref{subsubsec:SiPo} (where $M_c = M$) or one of the two sub-sampling strategies, the random single-point sub-sampling described in Section \ref{subsubsec:SiPo} or the patterned interrogation described in Section \ref{subsubsec:PatInter}.

\subsection{Image Reconstruction Strategies}  \label{subsec:ImageRecStrat}

Given the  compressed data $f^c$, there are in principle two strategies of how to reconstruct $p_{0}$: In \emph{two-step procedures}, we first reconstruct the detection plane pressure $\bar{p}$ from $f^c$ based on $f^c = C \bar{p}$ (\emph{data reconstruction}) and then use a standard PAT reconstruction for complete planar data (see, e.g. \cite{KuKu11}). In \emph{one-step procedures}, $p_{0}$ is reconstructed directly from $f^c$ using a model-based approach \eref{eq:Fwd}. While both approaches have advantages and disadvantages and novel two-step procedures have been introduced in \cite{HuZhBeArBeCo15,SaKrBeBuHa15,BeCoHuZhBeAr16}, the focus of this work is not to carry out a fair, detailed comparison between one-step and two-step approaches: We rather want to emphasize on the differences between simple, \emph{linear} reconstruction techniques and \emph{variational}, model-based reconstruction techniques, independent of the former being two-step and the latter being one-step procedures. To ease the following presentation, we first introduce the discrete PAT model we will use for numerical computations before discussing the details of the reconstruction techniques. 

\subsection{Discrete Forward Model}  \label{subsec:DiscFwdModel}

As all methods we examine directly rely on the wave equation \eref{eq:PATfwd}, we need a fast numerical method for 3D wave propagation with high spatial and temporal resolution. Our choice is the $k$-space pseudospectral time domain method \cite{MaSoLiTaNaWa01,CoKaArBe07,TrZhCo10} implemented in the \textit{k-Wave} Matlab Toolbox \cite{TrCo10}. For the following, it is only important that k-Wave discretizes $\Omega = [0,l_x] \times [0,l_y] \times [0,l_z]$ into $N = N_x \times N_y \times N_z$ voxels $\Delta_{x/y/z} = l_{x/y/z}/N_{x/y/z}$ and uses an explicit time stepping. The time step used will always be the same as used in the temporal sampling of the pressure field, i.e., $\Delta_t = \delta_t$ (cf. Section \ref{subsec:SubSam}). From now on, all variables used are discrete although we will continue to use the same notation for them. For instance, the discretization of the PAT forward operator is now a matrix $A \in \R^{N_y N_z M_t \times N}$, mapping the discrete initial pressure $p_0 \in \R^N$ to the pressure at the first layer of voxels in $x$ direction at the $M_t$ discrete time steps, as these voxels represent the detection plane. Note that we cannot construct $A$ explicitly, but rely on computing matrix-vector products with $A$ and $A^T$ using k-Wave. A detailed discussion of the implementation can be found in \cite{ArBeCoLuTr16}. The discrete sub-sampling operators $C$ map from the pressure-time series of the detection plane voxels to $f^c \in \R^{M_c M_t}$. The single point sampling operators (cf. Sections \ref{subsubsec:SiPo} and \ref{subsubsec:SiPo}) simply extract the pressure-time series at the sensor voxels (which are, in general, only a subset of the detection plane voxels). For the interrogation with binary sensing pattern constructed by the DMD (cf. Section \ref{subsec:CSFP}), \eref{eq:GenScan} can be implemented by multiplying a $M_c \times N_y N_z$ binary matrix $H$ with the vector of pressure values at the scanning voxels, separately for each time point $t$. Compressed sensing theory suggests that random Bernoulli matrices are optimal for many applications \cite{FoRa13}. As those are difficult to implement experimentally and computationally, we use scrambled Hadamard matrices which are known to have very similar properties compared to a random Bernoulli matrices and can be implemented very efficiently by a fast-fourier-transform-like operation \cite{DoLuNgTr12,FoRa13}. Note that as the entries of Bernoulli and Hadamard matrices take the values $\{-1,1\}$ and not $\{0,1\}$ as implemented by the DMD, we need to demean experimental data in a pre-processing step. Further details can be found in \cite{HuZhBeArBeCo15,BeCoHuZhBeAr16,HuZhBeArBeCo16}.

\subsection{Linear Back-Projection-Type Reconstructions}  \label{subsec:TwoStep}

As described above, we will use simple, linear two-step procedures to compare the more sophisticated variational methods to: First, we reconstruct the complete data by the pseudo-inverse of $C$: $C^\dagger f^c$ (cf. \eref{eq:Fwd}), where for the sub-sampling operators we consider here, $C^\dagger = C^T$. Then, we either multiply with $A^T \in \R^{N \times N_y N_z M_t}$ (called \termabb{back-projection}{BP} here), or with the discrete \termabb{time reversal}{TR} \cite{FiPa04,XuWa04b,JaLaOgTrCoZhJoPiPhMaLyPePuBe15} operator $A^{\triangleleft} \in \R^{N \times N_y N_z M_t}$:
\begin{eqnarray}
p_{\text{\tiny BP}} &\mydef A^T C^T f^c \label{eq:BP}\\
p_{\text{\tiny TR}} &\mydef A^{\triangleleft} C^T f^c \label{eq:TRnum}
\end{eqnarray}
The difference between  BP and TR (including enhanced variants of TR \cite{StUh09}), is discussed in more detail in \cite{ArBeCoLuTr16}. In summary, in the continuous setting, they differ in the way they introduce the time-reversed pressure time series at the detection plane: TR approaches use them as a time dependent Dirichlet boundary condition for the wave equation \eref{eq:PATfwd} while the adjoint approach introduces them as a time dependent source term without altering the boundary conditions.

\subsection{Variational Image Reconstruction}  \label{subsec:VarImRec}

\emph{Variational regularization} \cite{ScGrGrHaLe09} is a popular and well understood approach for approximating the solutions of ill-posed operator equations like \eref{eq:Fwd} in a reasonable and controlled manner: The regularized solution is defined as the minimizer of a suitable energy functional $\mathcal{E}$. Assuming that the additive noise, $\varepsilon$, is i.i.d. normally distributed, a reasonable approach is to solve 
\begin{equation}
p_{\lambda} \mydef \argminsub{p} \: \mathcal{E}(p) = \argminsub{p} \left\lbrace \frac{1}{2} \sqnorm{C A \, p - f^c} + \lambda \mathcal{J}(p) \right\rbrace, \label{eq:VarReg}
\end{equation}
to obtain a regularized solution $p_{\lambda}$. While the first term in the composite functional measures the misfit between measured and predicted data (\emph{data fidelity term}), $\mathcal{J}(p)$  has to render the minimization problem \eref{eq:VarReg} well-posed by ensuring existence, uniqueness and stability of $p_{\lambda}$ (\emph{regularization functional}). Furthermore, its choice can be used to penalize or constrain unwanted features of $p_{\lambda}$, thereby encoding \emph{a-priori} knowledge about the solution. The \emph{regularization parameter} $\lambda > 0$ controls the balance between both terms. \\
The first variational method we examine corresponds to classical, zeroth-order Tikhonov regularization, augmented by the physical constraint $p_{0} \geqslant 0$: 
\begin{equation}
p_{\text{\tiny L2+}} \mydef \argminsub{p \geqslant 0} \left\lbrace \frac{1}{2} \sqnorm{C A \, p - f^c} +  \lambda \sqnorm{p} \right\rbrace. \label{eq:L2+}
\end{equation}
As a second functional $\mathcal{J}(p)$, we examine the popular \termabb{total variation}{TV} energy, which is a discrete version of the \emph{total-variation} seminorm \cite{RuOsFa92,BuOs13}:
\begin{equation}
p_{\text{\tiny TV+}} \mydef \argminsub{p \geqslant 0} \left\lbrace \frac{1}{2} \sqnorm{C A \, p - f^c} +  \lambda \text{TV}(p) \right\rbrace. \label{eq:TV+}
\end{equation}
The energy $\text{TV}(p)$ measures the $\ell_1$ norm of the amplitude of the gradient field of $p$ (the details of its implementation are given in \ref{sec:TV}) and is a prominent example of non-smooth, \emph{edge-preserving} image reconstruction techniques, and, more generally, of \emph{spatial sparsity constraints}. While TV regularization has been used for PAT before (see, e.g., \cite{HuWaNiWaAn13}), our main interest in it arises from its results when applied to sub-sampled data for other imaging modalities
\cite{LuDoPa07,BrSaBu10,RiBeFlKa11,Mue13,HmKaKoLaNiSi13,BeGlHoScVa14,HaFaBeTaSeGlHo15,JoSi15}: TV regularization is often able to recover the object's main feature edges even for high sub-sampling factors. Therefore, we focus on this rather general regularization energy in this first PAT sub-sampling study and examine more specific functionals in future work. \\
As all of the involved functionals and constraints are convex, a variety of methods exist to solve \eref{eq:VarReg} computationally. In this work, we use an \emph{accelerated proximal gradient}-type method described in \ref{sec:Opti}.

\subsection{Bregman Iterations}  \label{subsec:BregIter}

A potential drawback of variational techniques like \eref{eq:VarReg} is that they inevitably lead to a systematic bias of the regularized solutions: The solution $p_{\lambda}$ moves from an un-biased data-fit towards a minimizer of $\mathcal{J}(p)$. Formally, let $\tilde{f}^c = C A p_{0}$ be the true, noise-free data and $\tilde{p}_{\lambda}$ the solution of \eref{eq:VarReg} for $f^c = \tilde{f}^c$. Then, by the minimizing properties of $\tilde{p}_{\lambda}$, we have 
\begin{equation}
\fl \quad \frac{1}{2} \sqnorm{C A \tilde{p}_{\lambda} - \tilde{f}^c} + \lambda \J(\tilde{p}_{\lambda}) \leqslant \frac{1}{2} \sqnorm{C A p_{0} - \tilde{f}^c} + \lambda \J(p_{0}) = \lambda \J(p_{0}),
\end{equation}
and thereby, $\J(\tilde{p}_{\lambda}) \leqslant \J(p_{0})$. For the TV energy, this bias manifests in the well-known, non-linear contrast loss of TV regularized solutions \cite{BuOs13}. For PAT, this systematic error poses a crucial limitation on the use of TV regularized PA images for quantitative analysis like QPAT studies. To overcome this drawback, an iterative enhancement of variational solutions by the \emph{Bregman iteration} \cite{Br67} was proposed in \cite{OsBuGoXuYi06}: For \eref{eq:VarReg}, they take the form
\begin{eqnarray}
p^{k+1}_{\lambda} &= \argminsub{p} \left\lbrace \frac{1}{2} \sqnorm{C A \, p - (f^c + b^k)} + \lambda \mathcal{J}(p) \right\rbrace, \label{eq:BregIterA}\\
b^{k+1}    &= b^k + \left( f^c - C A \, p^{k+1}_{\lambda} \right), \label{eq:BregIterB}
\end{eqnarray}
with $b^0 = 0$. This iteration has several attractive features \cite{BuOs13}: It solves the un-regularized problem
\begin{equation}
	\min_p \J(p)  \quad \text{subject to} \quad p \in \argminsub{q} \; \sqnorm{C A \, q - f^c}
\end{equation}
by solving a sequence of well-regularized problems \eref{eq:BregIterA} while the residual of iterates, $\norm{C A \, p^{k+1}_{\lambda} - f^c}$, is monotonically decreasing. The potential of Bregman iterations, in particular when used on sub-sampled data, has been demonstrated in \cite{BrSaBu10,Mue13,BeGlHoScVa14,HaFaBeTaSeGlHo15}. Note the difference between the use of the Bregman iteration in the \emph{Split Bregman method} \cite{GoOs09}, a method to solve problems like \eref{eq:VarReg} which is also known as the \emph{augmented Lagrangian method}, and the usage here which does not have an equivalent Lagrangian interpretation.


\section{Simulation Studies} \label{sec:SimStud}

We now examine the different inverse methods described in the previous section when applied to sub-sampled data from numerical phantoms.

\subsection{Realistic Numerical Phantom} \label{subsec:NumPhan}

\setlength{\fboxsep}{0pt}
\begin{figure}[tb]
\begin{minipage}[c]{0.725\textwidth}
\centering
\subfloat[][\label{subfig:NumPhanA}]{\includegraphics[width=\textwidth]{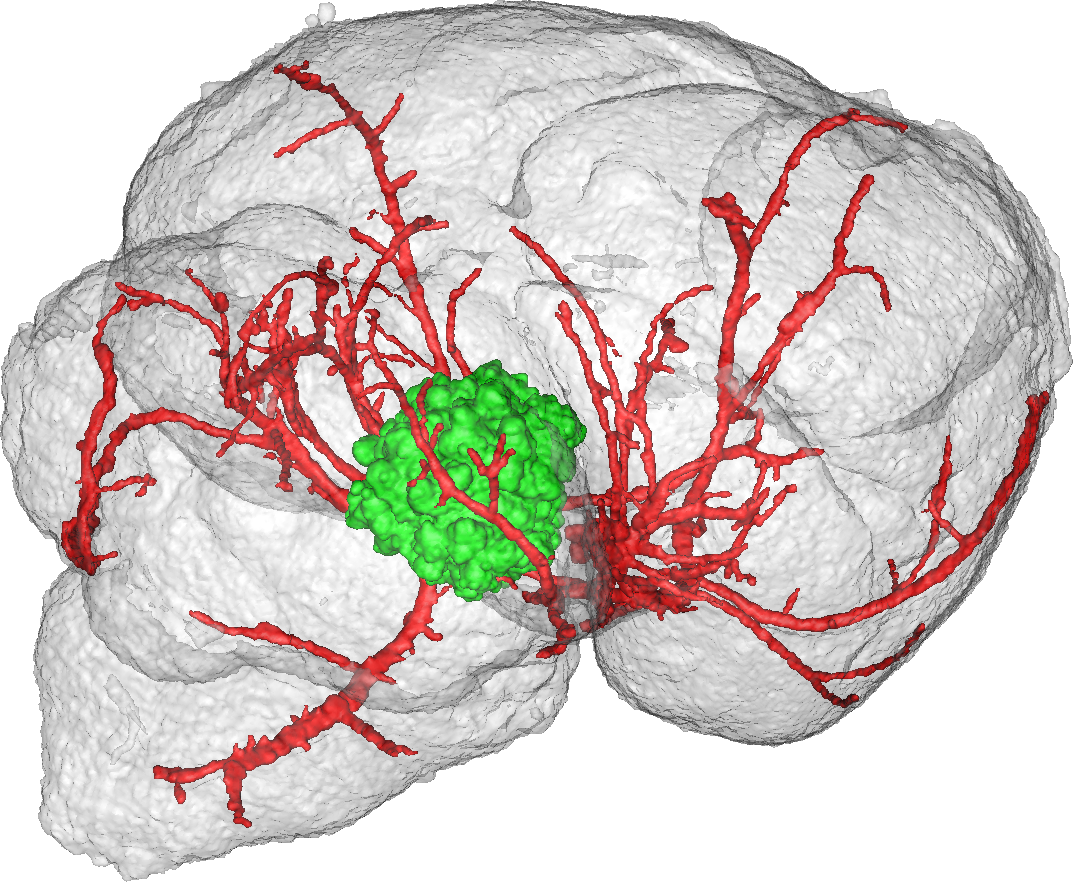}}
\end{minipage}
\hfill
\begin{minipage}[c]{0.225\textwidth}
\centering
\subfloat[][ \label{subfig:NumPhanB}]{\includegraphics[width=1\textwidth]{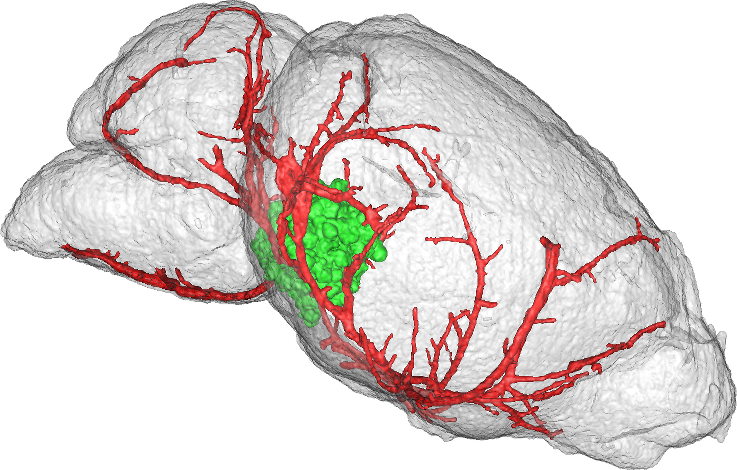}}\\
\subfloat[][ \label{subfig:NumPhanC}]{\includegraphics[width=1\textwidth]{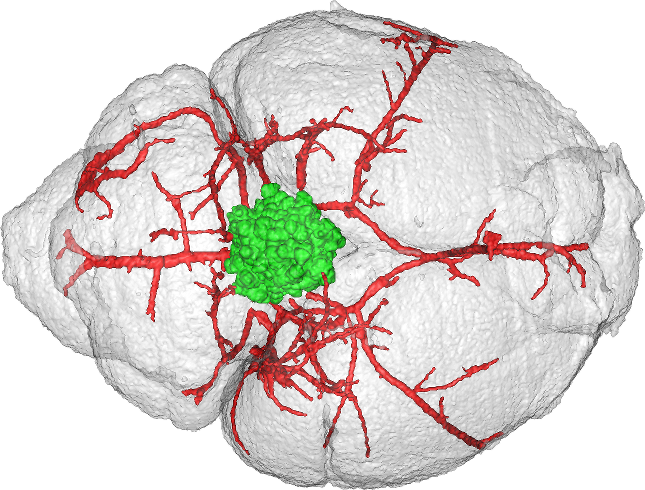}}\\
\subfloat[][ \label{subfig:NumPhanD}]{\includegraphics[width=1\textwidth]{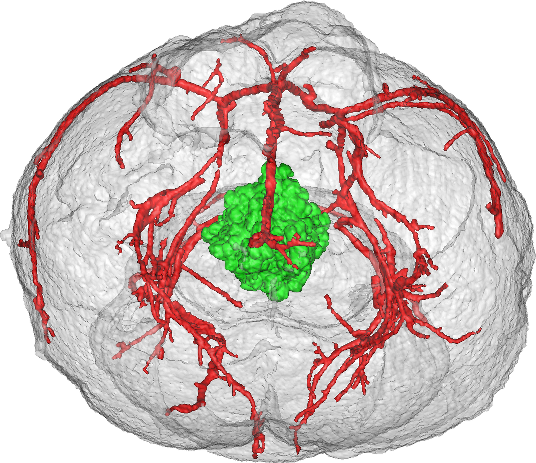}}
\end{minipage}
\caption{3D visualization (\protect\subref{subfig:NumPhanA}-\protect\subref{subfig:NumPhanD}: different views) of the segmentation used for the construction of the realistic numerical Phantom: Vasculature (red), artificial tumor tissue (green) and gray matter (gray). Volume rendering was carried out by \cite{SCIRun}.}
   \label{fig:NumPhan}
\end{figure}

\begin{figure}[tb]
   \centering
\subfloat[][X pro (y horiz, z vert)\label{subfig:VaTu1X}]{\includegraphics[width=0.32\textwidth]{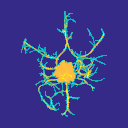}}
\hfill
\subfloat[][Y pro (x horiz, z vert)\label{subfig:VaTu1Y}]{\includegraphics[width=0.32\textwidth]{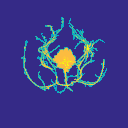}}
\hfill
\subfloat[][Z pro (x horiz, y vert)\label{subfig:VaTu1Z}]{\includegraphics[width=0.32\textwidth]{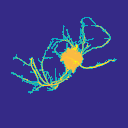}}
\\
\subfloat[][X pro (y horiz, z vert)\label{subfig:VaTu2X}]{\includegraphics[width=0.32\textwidth]{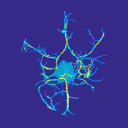}}
\hfill
\subfloat[][Y pro (x horiz, z vert)\label{subfig:VaTu2Y}]{\includegraphics[width=0.32\textwidth]{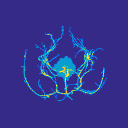}}
\hfill
\subfloat[][Z pro (x horiz, y vert)\label{subfig:VaTu2Z}]{\includegraphics[width=0.32\textwidth]{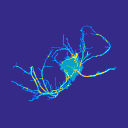}}
\caption{Maximum intensity projections of the realistic numerical phantoms (color map "parula", cf. Figure \ref{subfig:Parula}): \protect\subref{subfig:VaTu1X}-\protect\subref{subfig:VaTu1Z} Tumor1
\protect\subref{subfig:VaTu2X}-\protect\subref{subfig:VaTu2Z} Tumor2. In both cases, the detection plane corresponds to the top edge of the Y and Z projections.
}
   \label{fig:VaTu}
\end{figure}

While studies using numerical phantoms composed of simple geometrical objects can provide valuable insights to the basic properties of the inverse problem and reconstruction methods, it is often unclear how their results translate to experimental data from complex targets. In this section, we briefly describe the construction of a realistic numerical phantom that will be used in the main simulation studies. \\
The phantom is based on a segmentation of a micro-CT scan of a mouse brain ($533\times 400 \times 346$ voxel) into gray matter, vasculature and dura mater. The vasculature is morphologically closed (dilation followed
by erosion) using the 18-neighborhood as a convolution kernel. Thereafter, the vasculature is one connected component with respect to the 26-neighborhood. The whole segmentation is clipped to the bounding box of the vasculature leading to a size of $306\times 423 \times 345$. Next, a gray matter voxel in the central part of the segmentation is chosen uniformly at random. It is used as the seed point for the construction of an artificial cancer tissue inside the gray matter by a stochastic growth process that consists of an iterative application of morphological operations on the surface voxels. Figure \ref{fig:NumPhan} shows the final result of this construction. Note that vasculature and tumor tissue are non-intersecting. \\
For the studies in this work, the clipped volume is embedded into a cubic volume of $512^3$ voxels, centered in $y$ and $z$ direction and in two different heights in $x$ direction: In the first phantom, called \emph{Tumor1}, the distance between the detection plane $x = 0$ and the phantom is half of the distance between the plane $x = l_x$ and the phantom. In the second phantom, called \emph{Tumor2}, it is centered in $x$ direction, leading to a larger distance between sensor and target (cf. Figure \ref{fig:VaTu}). To construct $p_0$ from the segmentation, all vasculature voxels are given the value of $1$, whereas the tumor tissue is given the value of $0.7$ for Tumor1 and $0.3$ for Tumor2. Next, $p_0$ is down-sampled to the desired resolution $[N_x,N_y,N_z]$ by successive sub-averaging over $2\times 2 \times 2$ blocks. For Tumor1, we modify the resulting $p_0$ to obtain sharper boundaries: $p_0$ is normalized such that $\max(p_0) = 1$ and the intensity $p_{0,i}$ of all non-zero voxels $i$ is set to $(2 p_{0,i} + 1)/3$, thereby ensuring a contrast minimal value of $1/3$ between background and target. Figure \ref{fig:VaTu} shows \termabb{maximum intensity projections}{mxIP} of the resulting $p_0$.

\subsection{Simulation Studies with Tumor1} \label{subsec:SimVaTu1}

\begin{figure}[tb]
   \centering
\subfloat[][Tumor1, Y pro\label{subfig:EasyResPhan}]{\includegraphics[width=0.32\textwidth]{VascularTumorEasy128_parula-histCutOffPos_mIYpro.png}}
\hfill
\subfloat[][rSP-128\label{subfig:EasyResUni}]{\fbox{\includegraphics[width=0.32\textwidth]{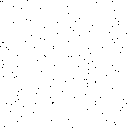}}}
\hfill
\subfloat[][one $\phi_j(y,z)$ from sHd-128 \label{subfig:EasyResPat}]{\fbox{\includegraphics[width=0.32\textwidth]{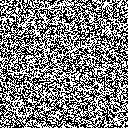}}}
\\[-5pt]
\subfloat[][TR, cnv., 30.99\si{\decibel}\label{subfig:EasyResTRf}]{\includegraphics[width=0.32\textwidth]{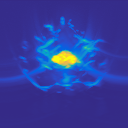}}
\hfill
\subfloat[][TR, rSP-128 29.47\si{\decibel}\label{subfig:EasyResTRUni}]{\includegraphics[width=0.32\textwidth]{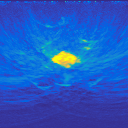}}
\hfill
\subfloat[][TR, sHd-128 15.44\si{\decibel}\label{subfig:EasyResTRPat}]{\includegraphics[width=0.32\textwidth]{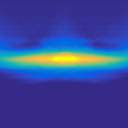}}
\\[-5pt]
\subfloat[][BP, cnv., 30.60\si{\decibel}\label{subfig:EasyResBPf}]{\includegraphics[width=0.32\textwidth]{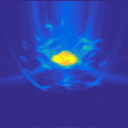}}
\hfill
\subfloat[][BP, rSP-128, 28.85\si{\decibel}\label{subfig:EasyResBPUni}]{\includegraphics[width=0.32\textwidth]{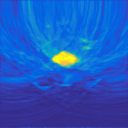}}
\hfill
\subfloat[][BP, sHd-128, 16.13\si{\decibel}\label{subfig:EasyResBPPat}]{\includegraphics[width=0.32\textwidth]{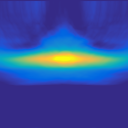}}
\caption{Tumor1 results (mxIP) of linear methods: \protect\subref{subfig:EasyResPhan} Phantom (cf. Figure \ref{fig:VaTu}), detection plane corresponds to the top edge. \protect\subref{subfig:EasyResUni} Visualization of the rSP-128 sub-sampling pattern: While each of the $M = 128 \times 128$ pixel corresponds to one possible scanning location, all black pixels correspond to one of the $M_c = 128$ random locations actually scanned. \protect\subref{subfig:EasyResPat} A $128 \times 128$ Hadamard matrix to symbolize the sHd-128 pattern (the actual matrix is of size $128 \times 16384$) \protect\subref{subfig:EasyResTRf}-\protect\subref{subfig:EasyResBPPat}: TR and BP results for conventional data (left column), rSP-128 (middle column) and sHd-128 (right column) and their corresponding PSNR in \si{\decibel}.}
   \label{fig:EasyLinRes}
\end{figure}

\begin{figure}[tb]
   \centering
\subfloat[][Tumor1, Y pro\label{subfig:EasyResPhan2}]{\includegraphics[width=0.32\textwidth]{VascularTumorEasy128_parula-histCutOffPos_mIYpro.png}}
\hfill
\subfloat[][rSP-128\label{subfig:EasyResUni2}]{\fbox{\includegraphics[width=0.32\textwidth]{VaTuE_128_siPi-uni-128-rs1.png}}}
\hfill
\subfloat[][one $\phi_j(y,z)$ from sHd-128 \label{subfig:EasyResPat2}]{\fbox{\includegraphics[width=0.32\textwidth]{scrambledHadamardSketch2.png}}}
\\[-5pt]
\subfloat[][L2+, cnv., 33.44\si{\decibel}\label{subfig:EasyResL2f}]{\includegraphics[width=0.32\textwidth]{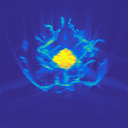}}
\hfill
\subfloat[][L2+, rSP-128, 31.25\si{\decibel}\label{subfig:EasyResL2Uni}]{\includegraphics[width=0.32\textwidth]{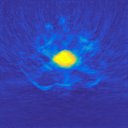}}
\hfill
\subfloat[][L2+, sHd-128, 33.54\si{\decibel}\label{subfig:EasyResL2Pat}]{\includegraphics[width=0.32\textwidth]{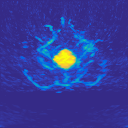}}
\\[-5pt]
\subfloat[][TV+, cnv., 33.84\si{\decibel}\label{subfig:EasyResTVf}]{\includegraphics[width=0.32\textwidth]{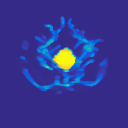}}
\hfill
\subfloat[][TV+, rSP-128, 33.49\si{\decibel}\label{subfig:EasyResTVUni}]{\includegraphics[width=0.32\textwidth]{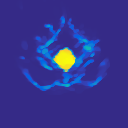}}
\hfill
\subfloat[][TV+, sHd-128, 33.64\si{\decibel}\label{subfig:EasyResTVPat}]{\includegraphics[width=0.32\textwidth]{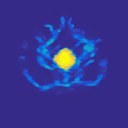}}
\\[-5pt]
\subfloat[][TV+Br, cnv., 34.06\si{\decibel}\label{subfig:EasyResTVBrf}]{\includegraphics[width=0.32\textwidth]{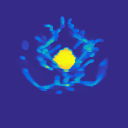}}
\hfill
\subfloat[][TV+Br, rSP-128, 33.39\si{\decibel}\label{subfig:EasyResTVBrUni}]{\includegraphics[width=0.32\textwidth]{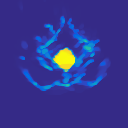}}
\hfill
\subfloat[][TV+Br, sHd-128, 33.41\si{\decibel}\label{subfig:EasyResTVBrPat}]{\includegraphics[width=0.32\textwidth]{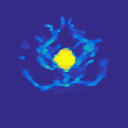}}
\caption{Tumor1 results (mxIP) of variational methods: \protect\subref{subfig:EasyResPhan2}-\protect\subref{subfig:EasyResPat2} cf. Figure \ref{fig:EasyLinRes} \protect\subref{subfig:EasyResL2f}-\protect\subref{subfig:EasyResTVBrPat}: L2+, TV+ and TV+Br results for conventional data (left column), rSP-128 (middle column) and sHd-128 (right column) and their corresponding PSNR in \si{\decibel}.}
   \label{fig:EasyVarRes}
\end{figure}

We first examine the inverse reconstructions using Tumor1 with $N = 128^3$ for "inverse crime" data \cite{KaSo07}, which means that we assume that we have exact knowledge of all physical parameters, which are summarized in Table \ref{tbl:ModelPara}, and use the same model for both data simulation and reconstruction. In addition, we sampled the data with the spatial spacing given by the Nyquist criterion: The spatial sampling intervals $\delta_{y/z}$ coincide with the spatial spacing of the computational grid $\Delta_{y/z}$ and are fine enough to capture all relevant spatial frequencies of the incident photoacoustic field (cf. Section \ref{subsec:BackNyquist}): $p_0$ was pre-smoothed to ensure that its discrete approximation by the truncated Fourier basis used in k-Wave was non-oscillatory. This is reflected in a sharp drop of the power spectrum of the pressure time series around 4.8 \si{\mega \hertz}, which (cf. Section \ref{subsec:BackNyquist}) corresponds to a spatial Nyquist rate equal to the spatial spacing $\delta_{y/z} = \Delta_{y/z} = 156.25$ \si{\micro \meter}. White noise with a standard deviation of $\sigma = 0.001$ is added to the clean data $C A p_0$, leading to a \termabb{signal-to-noise ratio}{SNR} of $18.63$ \si{\decibel}.\\
For all computed solutions $p$, voxels with negative pressure and all sensor voxels have been set to $0$ in a post-processing step. We will mainly rely on a visual comparison of the results via maximum intensity projections along the $y$ direction (cf. Figures \ref{subfig:VaTu1Y}, \ref{subfig:VaTu2Y}). Unless stated otherwise, the color scale of each figure is determined independently. It ranges from $0$ to the value separating the $100/256\% \approx 0.39\%$ largest values of $p$ from the smaller values. This clipping is necessary to avoid that a few large outliers determine the contrast of the image and complicate the comparison between different methods. In addition to the visual comparison, we report the \termabb{mean-squared-error}{MSE} between the reconstructed solution $p$ and $p_0$, i.e., $\|p-p_0\|_2^2/N$, in the conventional logarithmic scaling termed \termabb{peak-signal-to-noise ratio}{PSNR}:
\begin{equation}
\text{PSNR}(p,p_0) \mydef - 10 \log_{10} \left( \frac{1}{N} \|p-p_0\|_2^2 \right) \label{eq:PSNR} \\ 
\end{equation}
$p$ might have a different scaling compared to $p_0$ and contain small scale noise. As this typically does not influence the evaluation of a human observer, we also do not want to account for it when computing the PSNR. Therefore, we first rescale and threshold $p$ and $p_0$,
\begin{equation*}
\fl \qquad \tilde{p}  = \text{thres}\left( \frac{p}{\|p\|_\infty},0.1\right), \nonumber
 \; \tilde{p_0} \;  \text{accordingly}, \quad \text{where} \quad 
 \text{thres}(v,\alpha) = \cases{v \quad \text{if} \quad v \geqslant \alpha \\
0 \quad \text{else}
},
\end{equation*}
and then compute $\text{PSNR}(\tilde{p},\tilde{p_0})$ by \eref{eq:PSNR}.\\
Figure \ref{fig:EasyLinRes} shows the results of TR and BP for using rSP-128 or sHd-128 sub-sampling, i.e., accelerated by a factor of $M_{sub} = 128$ compared to the conventional scan. This high sub-sampling factor leads to unsatisfactory reconstructions. The different sub-sampling strategies lead to a different appearance of back-propagated noise and sub-sampling artifacts in the reconstructed images: In the rSP-128 case, both features and noise are back-propagated along the surfaces of spheres centered on the scanning locations while in the sHd-128 case, the non-local nature of the patterned interrogation leads to back-propagation along planes parallel to the detection plane. Figure \ref{fig:EasyVarRes} shows the corresponding results of the classical Tikhonov regularization \eref{eq:L2+}, denoted by "L2+", TV regularization \eref{eq:TV+}, denoted by "TV+" and of the Bregman iteration \eref{eq:BregIterA}-\eref{eq:BregIterB} applied to  TV regularization, denoted by "TV+Br". Despite the high sub-sampling factor, TV+ and TV+Br are still able to reconstruct the main structures of the phantom without excessive image noise. \\
In the results shown, the regularization parameter $\lambda$  for L2+ and TV+ was chosen by the \termabb{discrepancy principle}{DP}: For the general variational problem \eref{eq:VarReg}, the discrepancy principle selects $\lambda$ such that 
\begin{equation}
\frac{\norm{C A \, p_{\lambda} - f^c}}{\sqrt{M_c} \sigma}  = \kappa,
 \label{eq:DiscPr}
\end{equation}
for $\kappa \geqslant 1$. The DP is based on the heuristic argument, that the regularized solution $p_{\lambda}$ should explain the data $f^c$ no more than up to the noise level, which is assumed to be known. As the residual is monotonically increasing in $\lambda$, the DP is robust and easy-to-implement. We chose a simple interval-based method that linearly interpolates the left hand side of \eref{eq:DiscPr} (the \emph{discrepancy} of the data) in the current search interval. It was terminated when $\kappa = 1.25$ was reached within a tolerance of $0.01$. The Bregman iterations were started with $\lambda_{\text{\tiny TV+Br}}$ = $10 \lambda_{\text{\tiny TV+}}$, where $\lambda_{\text{\tiny TV+}}$ is the regularization parameter found for TV+, and stopped as soon as the discrepancy of $p^{k+1}_{\lambda}$ falls below $\kappa$ (recall that the residual, and thereby the discrepancy monotonically decreases with $k$, cf. Section \ref{subsec:BregIter}). This typically happens after about 10 Bregman iterations.

\subsubsection{Enhancement through Bregman Iterations} \label{subsubsec:EnhanceBregIter}

\begin{figure}[tb]
   \centering
\subfloat[][TV+, cnv. data \label{subfig:BeneBrTV}]{\includegraphics[width=0.24\textwidth]{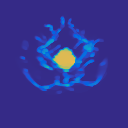}}
\hfill
\subfloat[][TV+Br, cnv. data\label{subfig:BeneBrTVBr}]{\includegraphics[width=0.24\textwidth]{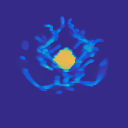}}
\hfill
\subfloat[][$(p_{\text{\tiny TV+Br}}-p_{\text{\tiny TV+}})_+$, cnv. data\label{subfig:BeneBrDiffPos}]{\fbox{\includegraphics[width=0.24\textwidth]{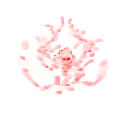}}}
\hfill
\subfloat[][$(p_{\text{\tiny TV+Br}}-p_{\text{\tiny TV+}})_-$, cnv. data\label{subfig:BeneBrDiffNeg}]{\fbox{\includegraphics[width=0.24\textwidth]{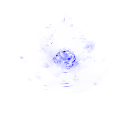}}}
\\
\subfloat[][TV+, rSP-128 \label{subfig:BeneBrTVUni}]{\includegraphics[width=0.24\textwidth]{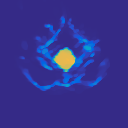}}
\hfill
\subfloat[][TV+Br, rSP-128\label{subfig:BeneBrTVBrUni}]{\includegraphics[width=0.24\textwidth]{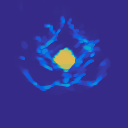}}
\hfill
\subfloat[][$(p_{\text{\tiny TV+Br}}-p_{\text{\tiny TV+}})_+$, rSP-128\label{subfig:BeneBrDiffPosUni}]{\fbox{\includegraphics[width=0.24\textwidth]{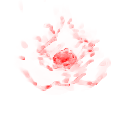}}}
\hfill
\subfloat[][$(p_{\text{\tiny TV+Br}}-p_{\text{\tiny TV+}})_-$, rSP-128\label{subfig:BeneBrDiffNegUni}]{\fbox{\includegraphics[width=0.24\textwidth]{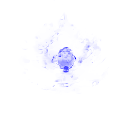}}}
\caption{Contrast comparison between TV+ and TV+Br solutions (cf. Figure \ref{fig:EasyVarRes}): The images in the left two columns share the same color scale. The third and forth column show mxIPs of the positive (red scale) and negative (blue scale) part of the difference $p_{\text{\tiny TV+Br}}-p_{\text{\tiny TV+}}$. 
}
\label{fig:BeneBr}
\end{figure}

The Bregman iteration was introduced to compensate for the systematic contrast loss of TV regularized solutions. Figure \ref{fig:BeneBr} compares TV+ and TV+Br solutions using the same color scaling and additionally shows mxIPs of the positive and negative parts of the difference between TV+Br and TV+. The difference plots demonstrate that using Bregman iterations especially improves the contrast of the small scale vessel structures and that the benefit is more pronounced for sub-sampled data compared to conventional data.

\subsection{Simulation Studies with Tumor2} \label{subsec:SimVaTu2}

The good quality of, e.g., the TV+Br results for the very high sub-sampling factor of $128$ have to be interpreted with care: The phantom Tumor1 used was intentionally easy to reconstruct as it was close to the sensors and had high contrast (cf. Figure \ref{fig:VaTu}). In addition, the data was created by the same forward model used in the reconstruction, which is known as committing an "inverse crime" \cite{KaSo07}, and the conventional data was sampled finely enough to fulfill the Nyquist criterion. The phantom Tumor2 was designed to carry out simulation studies that more accurately reproduce the challenges of experimental data scenarios. 

\subsubsection{Inverse Crimes} \label{subsubsec:InvCrimes}

\begin{figure}[tb]
   \centering
\subfloat[][$c_0 + \tilde{c}$\label{subfig:SoundSpeed}]{\includegraphics[height=0.4\textwidth]{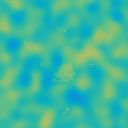}}
\hskip 10pt
\subfloat[][\label{subfig:Parula}]{\fbox{\includegraphics[height=0.4\textwidth]{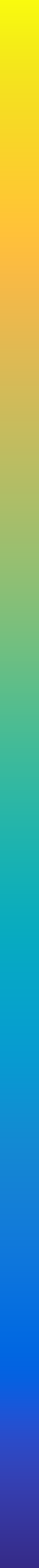}}}
\hfill
\subfloat[][pressure-time series\label{subfig:PressureTimeCourse}]{\includegraphics[height=0.4\textwidth]{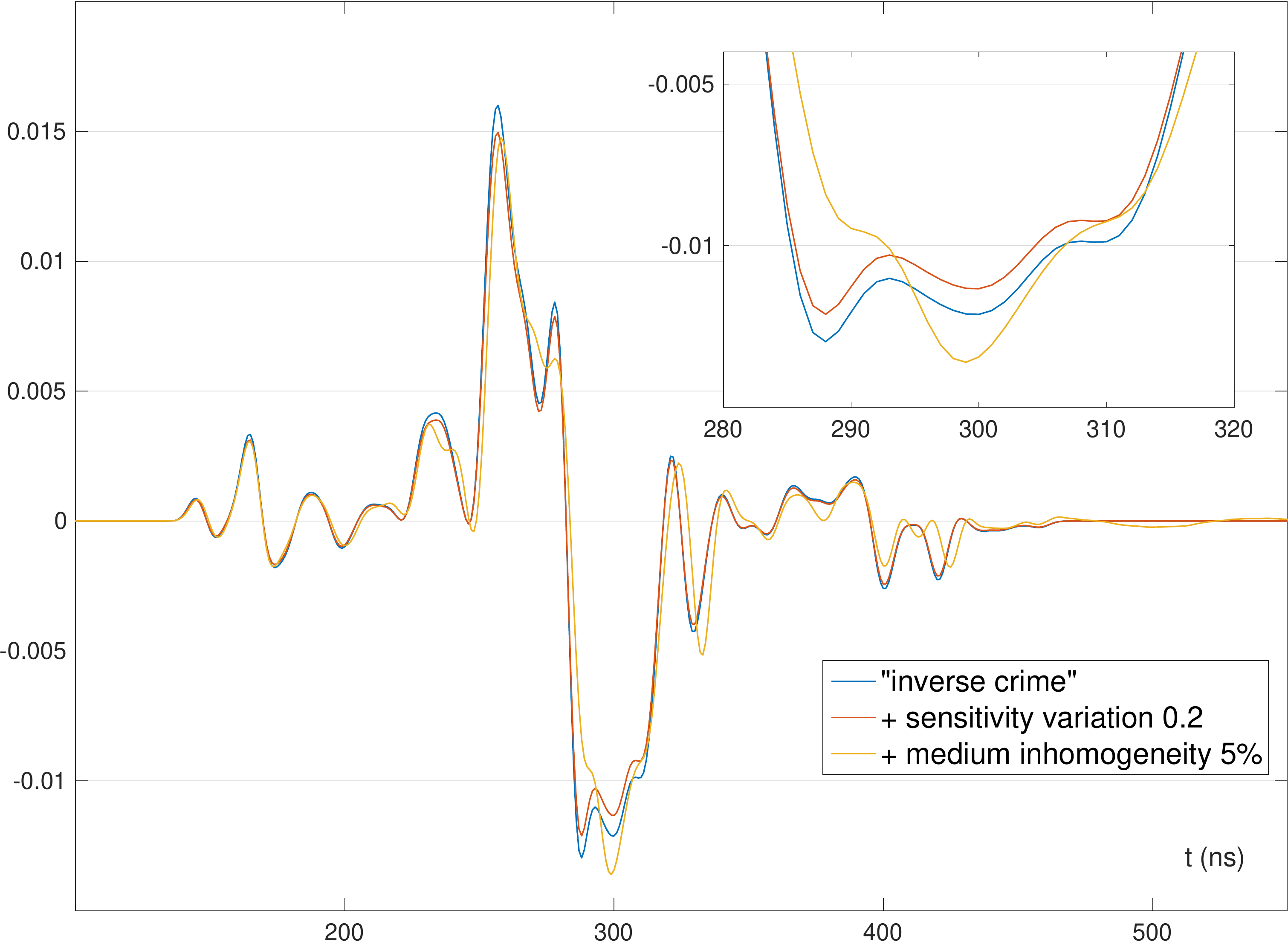}}
\caption{
\protect\subref{subfig:SoundSpeed} $x$-slice of the sound speed used to generate "no-inverse-crime" data for the Tumor2 phantom \protect\subref{subfig:Parula} color scale, range: 1350-1650 \si{\meter\per\second}, \protect\subref{subfig:PressureTimeCourse} noise-free pressure-time series to demonstrate the effect of sensitivity and sound speed variation (inset zooms into particular section of the plot).
}
   \label{fig:IC}
\end{figure}

While inverse crimes are, to a certain extend, unavoidable when carrying out simulation studies, they make it more difficult to extrapolate the results obtained to experimental data. In this section, we discuss how to bridge this gap by modifying the model used for the data generation:
\begin{itemize}
\item The sensitivity of the FP sensor varies spatially \cite{ZhLaBe08}. While this effect can be mitigated by calibration and data pre-processing procedures, a residual uncertainty remains that we will model by modifying the discrete (static) data generation model to 
\begin{equation}
f^c = C W_{s} A p_{0} + \varepsilon \label{eq:FwdSenVar}
\end{equation}
where $W_{s}$ is a diagonal matrix that multiplies the pressure-time course of each voxel $i$ in the $x=1$ plane by a random variable $w_i$ following a centered log-normal distribution:
\begin{equation}
w_i = \exp(\sigma_{s} X_i), \qquad X_i \sim \mathcal{N}(0,1) \label{eq:LogNorm}
\end{equation}
We choose $\sigma_{s} = 0.2$ (90\% of $w_i$ are in $[0.72,1.39]$).
\item In a similar spirit, we assume that we only have a rough estimate of the statistics of $\varepsilon$, e.g., from baseline measurements, and can, therefore, only approximately decorrelate the data before the inversion. The residual uncertainty about the noise variance per sensor voxel is modeled by replacing the additional noise term $\varepsilon$ by $W_{n} \varepsilon$, where $W_{n}$ is constructed like $W_s$ with $\sigma_n = 0.1$ (90\% of $w_i$ are in $[0.85,1.18]$). Keeping $\sigma = 0.001$, as the standard deviation of $\varepsilon$, we end up with an average SNR of $9.51$ (the value is considerable lower than for Tumor1 due to the larger distance to the sensors and the lower contrast of Tumor2).
\item Although we assume that the medium we image is sufficiently homogeneous to assume a constant sound speed $c_0$ in the inverse reconstruction, the real sound speed will slightly vary, especially between different tissue types. We use $c_0 + \tilde{c}$ for the data generation, where $\tilde{c}$ is constructed by adding a smooth, normalized Gaussian random field and the normalized initial pressure $p_0$ (as it represents a tissue different from the background). Then, $\tilde{c}$ is centered and scaled such that its mean is $0$ and its maximal absolute value is $0.05 c_0$ (a sound speed variation of $5$ \% is not unusual for soft tissue \cite{JiWa06}). The resulting sound speed is shown in Figure \ref{subfig:SoundSpeed}.
\item Among the many other ways to modify the data generation model are to include acoustic absorption, inhomogeneous illumination, acoustic reflections, baseline shifts and drifts in the pressure-time series, correlated noise and corrupted channels. We leave these extensions to future studies.  
\end{itemize}
Figure \ref{subfig:PressureTimeCourse} compares the noise-free pressure-time series of a single voxel after adding sensitivity and sound speed variations (adding noise and noise variation complicates a visual comparison).

\subsubsection{Nyquist criterion} \label{subsubsec:NyquistVaTu2}

In contrast to the previous studies, we now compare to conventional data acquired on a regular grid having a grid spacing corresponding to twice the length determined by the Nyquist criterion (cf. Section \ref{subsec:SimVaTu1}): $\delta_{y/z} = 2 \Delta_{y/z} = 312.5$ \si{\micro \meter}, i.e., we model the conventional data by extracting the pressure-time series at a sub-set of the $128^2$ detection plane voxels forming a regular grid with spacing 2. This is closer to the measured datasets we will examine in Section \ref{sec:ExpData}. All acceleration factors are defined with respect to the total number of locations of the regular grid which is given by $M = 128^2/2^2 = 4096$. Note that it now also makes sense to consider rSP-1 and sHd-1 as "sub-sampling" pattern: Although $M_{c} = M$, the data they measure cannot be converted to the conventionally sampled data, as was the case with Tumor1.

\subsubsection{Reconstruction Results} \label{subsubsec:ResVaTu2}
 
\begin{figure}[tb]
   \centering
\subfloat[][Tumor2, Y mxIP\label{subfig:noICResPhan}]{\includegraphics[width=0.32\textwidth]{VascularTumor128_parula-histCutOffPos_mIYpro.png}}
\hfill
\subfloat[][rSP-16\label{subfig:noICResUniPat}]{\fbox{\includegraphics[width=0.32\textwidth]{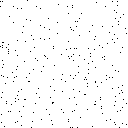}}}
\hfill
\subfloat[][one $\phi_j(y,z)$ from sHd-16\label{subfig:noICResPat}]{\fbox{\includegraphics[width=0.32\textwidth]{scrambledHadamardSketch2.png}}}
\\[-5pt]
\subfloat[][L2+, cnv., 33.46\si{\decibel}\label{subfig:noICL2}]{\includegraphics[width=0.32\textwidth]{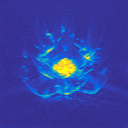}}
\hfill
\subfloat[][L2+, rSP-16, 31.20\si{\decibel}\label{subfig:noICL2Uni}]{\includegraphics[width=0.32\textwidth]{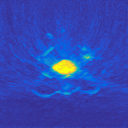}}
\hfill
\subfloat[][L2+, sHd-16, 33.57\si{\decibel}\label{subfig:noICL2Pat}]{\includegraphics[width=0.32\textwidth]{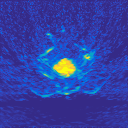}}
\\[-5pt]
\subfloat[][TV+, cnv., 35.11\si{\decibel}\label{subfig:noICTV}]{\includegraphics[width=0.32\textwidth]{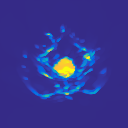}}
\hfill
\subfloat[][TV+, rSP-16, 33.06\si{\decibel}\label{subfig:noICTVUni}]{\includegraphics[width=0.32\textwidth]{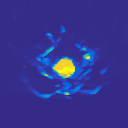}}
\hfill
\subfloat[][TV+, sHd-16, 34.59\si{\decibel}\label{subfig:noICTVPat}]{\includegraphics[width=0.32\textwidth]{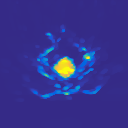}}
\\[-5pt]
\subfloat[][TV+Br, cnv., 35.78\si{\decibel}\label{subfig:noICBr}]{\includegraphics[width=0.32\textwidth]{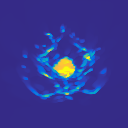}}
\hfill
\subfloat[][TV+Br, rSP-16, 34.16\si{\decibel}\label{subfig:noICBrUni}]{\includegraphics[width=0.32\textwidth]{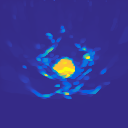}}
\hfill
\subfloat[][TV+Br, sHd-16, 34.27\si{\decibel}\label{subfig:noICBrPat}]{\includegraphics[width=0.32\textwidth]{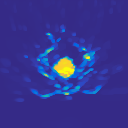}}
\caption{Tumor2 results (mxIP) of variational methods: \protect\subref{subfig:noICResPhan}-\protect\subref{subfig:noICResPat} cf. Figure \ref{fig:EasyLinRes} \protect\subref{subfig:noICL2}-\protect\subref{subfig:noICBrPat}: L2+, TV+ and TV+Br results for conventional data (left column), rSP-16 (middle column) and sHd-16 (right column) and their corresponding PSNR in \si{\decibel}.}
   \label{fig:noICAllRes}
\end{figure}

\begin{figure}[t]
   \centering
\subfloat[][TV+Br, cnv., 35.78\si{\decibel}\label{subfig:noICBrFull}]{\includegraphics[width=0.28\textwidth]{{VaTu_128_siPi-grid-4096_NoIC-1-0.2-0.05-0.1_std1.0e-03RS1_TV+Br_lam3.31e-03-0disc_PrxGr-mxIt100_parula-histCutOffPos_mIY}.png}}
\hfill
\subfloat[][rSP-1, 35.24\si{\decibel}\label{subfig:noICBr1Uni}]{\includegraphics[width=0.28\textwidth]{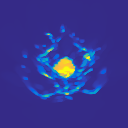}}
\hfill
\subfloat[][sHd-1, 35.97\si{\decibel}\label{subfig:noICBr1Pat}]{\includegraphics[width=0.28\textwidth]{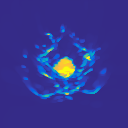}}
\\[-7pt]
\captionsetup[subfloat]{labelformat=empty}
\subfloat[][]{\includegraphics[width=0.28\textwidth]{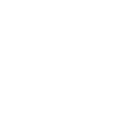}}
\captionsetup[subfloat]{labelformat=parens}
\hfill
\subfloat[][rSP-4, 35.14\si{\decibel}\label{subfig:noICBr4Uni}]{\includegraphics[width=0.28\textwidth]{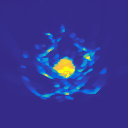}}
\hfill
\subfloat[][sHd-4, 35.32\si{\decibel}\label{subfig:noICBr4Pat}]{\includegraphics[width=0.28\textwidth]{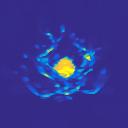}}
\\[-7pt]
\captionsetup[subfloat]{labelformat=empty}
\subfloat[][]{\includegraphics[width=0.28\textwidth]{dummy.png}}
\captionsetup[subfloat]{labelformat=parens}
\hfill
\subfloat[][rSP-8, 34.19\si{\decibel}\label{subfig:noICBr8Uni}]{\includegraphics[width=0.28\textwidth]{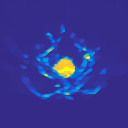}}
\hfill
\subfloat[][sHd-8, 34.45\si{\decibel}\label{subfig:noICBr8Pat}]{\includegraphics[width=0.28\textwidth]{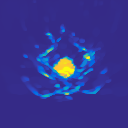}}
\\[-7pt]
\captionsetup[subfloat]{labelformat=empty}
\subfloat[][]{\includegraphics[width=0.28\textwidth]{dummy.png}}
\captionsetup[subfloat]{labelformat=parens}
\hfill
\subfloat[][rSP-16, 34.16\si{\decibel}\label{subfig:noICBr16Uni}]{\includegraphics[width=0.28\textwidth]{{VaTu_128_siPi-uni-256-rs1_NoIC-1-0.2-0.05-0.1_std1.0e-03RS1_TV+Br_lam4.82e-04-0disc_PrxGr-mxIt100_parula-histCutOffPos_mIY}.png}}
\hfill
\subfloat[][sHd-16, 34.27\si{\decibel}\label{subfig:noICBr16Pat}]{\includegraphics[width=0.28\textwidth]{{VaTu_128_fH256-rs1_NoIC-1-0.2-0.05-0.1_std1.0e-03RS1_TV+Br_lam4.41e-04-0disc_PrxGr-mxIt100_parula-histCutOffPos_mIY}.png}}
\\[-7pt]
\captionsetup[subfloat]{labelformat=empty}
\subfloat[][]{\includegraphics[width=0.28\textwidth]{dummy.png}}
\captionsetup[subfloat]{labelformat=parens}
\hfill
\subfloat[][rSP-32, 33.38\si{\decibel}\label{subfig:noICBr32Uni}]{\includegraphics[width=0.28\textwidth]{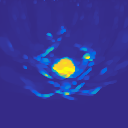}}
\hfill
\subfloat[][sHd-32, 33.65\si{\decibel}\label{subfig:noICBr32Pat}]{\includegraphics[width=0.28\textwidth]{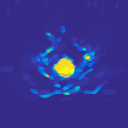}}
\caption{TV+Br results for Tumor2 phantom, comparison between different sub-sampling factors $M_{sub}$ }
   \label{fig:noICAllTVBrRes}
\end{figure}

As the results of TR and BP are of similar quality as for Tumor1 (cf. Figure \ref{fig:EasyLinRes}) we omit them here, and concentrate on the results of the variational methods. We again used the DP ($\kappa = 1.25$) to select the regularization parameter. Figure \ref{fig:noICAllRes} compares them for rSP-16 and sHd-16 sub-sampling: TV+Br, again, leads to the best results by visual impression and PSNR. In Figure \ref{fig:noICAllTVBrRes}, we therefore examine the influence of $M_{sub}$ on the reconstructed images only for TV+Br: Up to $M_{sub} = 16$, the rSP-based TV+Br reconstructions only slightly deteriorate. From $M_{sub} = 16$ to $M_{sub} = 32$, however, a clear degradation is visible. For the sHd sub-sampling, the image quality remains acceptable up to $M_{sub} = 32$.

\subsubsection{Influence of the Spatial Sub-Sampling Pattern}

\begin{figure}[tb]
   \centering
\subfloat[][Tumor2, Y mxIP\label{subfig:UniVsGridPhan}]{\includegraphics[width=0.32\textwidth]{VascularTumor128_parula-histCutOffPos_mIYpro.png}}
\hfill
\subfloat[][rSP-16\label{subfig:UniVsGridUniSam}]{\fbox{\includegraphics[width=0.32\textwidth]{VaTu_128_siPi-uni-256-rs1.png}}}
\hfill
\subfloat[][gSP-16\label{subfig:UniVsGridGrSam}]{\fbox{\includegraphics[width=0.32\textwidth]{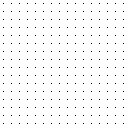}}}
\\[-5pt]
\subfloat[][TV+Br, cnv., 35.78\si{\decibel}\label{subfig:UniVsGridFulll}]{\includegraphics[width=0.32\textwidth]{{VaTu_128_siPi-grid-4096_NoIC-1-0.2-0.05-0.1_std1.0e-03RS1_TV+Br_lam3.31e-03-0disc_PrxGr-mxIt100_parula-histCutOffPos_mIY}.png}}
\hfill
\subfloat[][rSP-16, 34.16\si{\decibel}\label{subfig:UniVsGridUni}]{\includegraphics[width=0.32\textwidth]{{VaTu_128_siPi-uni-256-rs1_NoIC-1-0.2-0.05-0.1_std1.0e-03RS1_TV+Br_lam4.82e-04-0disc_PrxGr-mxIt100_parula-histCutOffPos_mIY}.png}}
\hfill
\subfloat[][gSP-16, 33.92\si{\decibel}\label{subfig:UniVsGridGr}]{\includegraphics[width=0.32\textwidth]{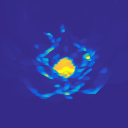}}
\caption{Influence of the spatial sub-sampling pattern in single point sub-sampling: \protect\subref{subfig:UniVsGridPhan} Phantom \protect\subref{subfig:UniVsGridUniSam}-\protect\subref{subfig:UniVsGridGrSam} visualizations of rSP-16 and gSP-16 \protect\subref{subfig:UniVsGridFulll}-\protect\subref{subfig:UniVsGridGr} TV+Br reconstructions.}
   \label{fig:UniVsGrid}
\end{figure}

As described in Section \ref{subsubsec:SiPo}, a random partition of the scanning locations was chosen as the main single point sub-sampling pattern to avoid unintended systematic artifacts like aliasing and was furthermore inspired by several results in compressed sensing theory \cite{FoRa13}. In Figure \ref{fig:UniVsGrid}, we compare this random pattern to using a regular sub-sampling pattern based on a coarse grid, which we denote as gSP-$M_{sub}$. The results show that the concrete choice of the single point sub-sampling pattern seems to be a minor influence, compared to, e.g., the choice of the inverse method used. We leave a more detailed examination and the design of optimal (dynamic) sampling pattern for further research. 

\section{Experimental Data} \label{sec:ExpData}

In this section, we examine the sub-sampling strategies for three example experimental data sets. To have a ground truth, sub-sampling was only carried out artificially: For each experiment, a conventional data set was acquired first, and sub-sampled data sets were produced from this data thereafter.

\subsection{Single Point Sub-Sampling - Dynamic Phantom} \label{subsec:Knot}

\begin{figure}[t]
   \centering
\subfloat[][Dynamic phantom "Knot" \label{subfig:ThKnSetup}]{\includegraphics[height=0.35\textwidth]{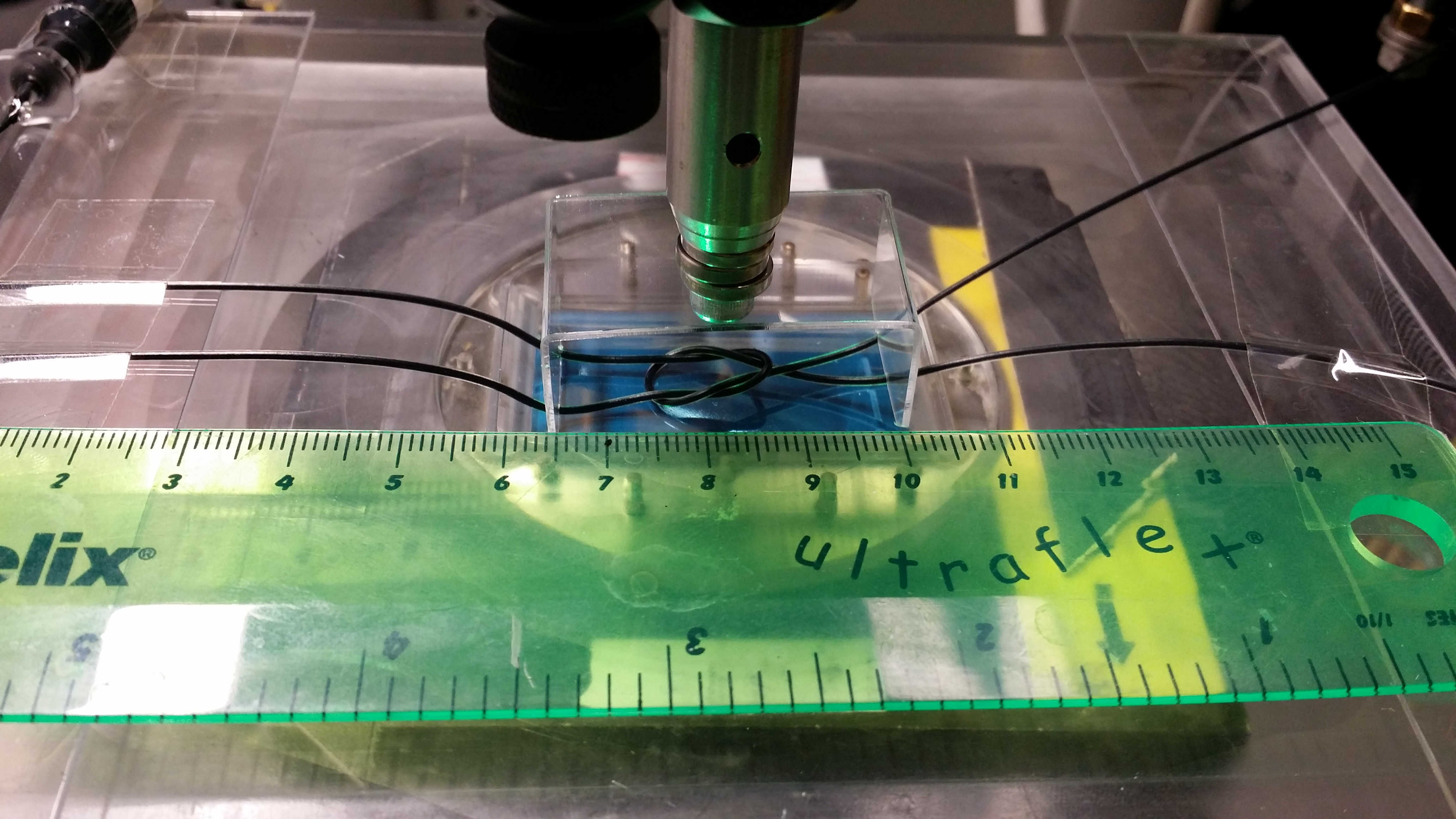}}
\hfill
\subfloat[][Static phantom "Hair" \label{subfig:HrKnSetup}]{\includegraphics[height=0.35\textwidth]{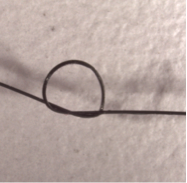}}
\caption{Experimental phantoms.}
   \label{fig:ExpPhantoms}
\end{figure}

We start with data acquired by a conventional single point FP scanner, measuring a pseudo-dynamic experimental phantom, which we call "Knot". 

\subsubsection{Setup} \label{subsubsec:KnotSetup}

Figure \ref{subfig:ThKnSetup} shows the experimental setup: Two polythene tubes were filled with 10\% and 100\% ink and interleaved to form a knot. One of the loose ends was tied to a motor shaft (top right corner of Figure \ref{subfig:ThKnSetup}) while the other three ends were fixated. The pseudo-dynamic data was acquired in a stop-motion style: A conventional scan (duration $\sim$15\si{\minute}) was performed while the whole arrangement was fixed. Then, the motor shaft was rotated by a stationary angle, causing the knot to tighten and to move into the direction of the motor, and a new scan was performed. This way, a conventional data set comprising $45$ frames was acquired. The tubes were immersed in a $1$\% Intralipid solution with de-ionised water. The excitation laser pulses used had a wavelength of $1064$\si{\nano\meter}, an energy of around $20$\si{\milli\joule} and were delivered a rate of $20$\si{\hertz}. A conventional scan consisted of $134 \times 133$ locations ($\delta_x = \delta_y = 150$\si{\micro\meter}), measured over $M_t = 625$ time points with a resolution of $\delta_t = 12$\si{\nano\second}.

\subsubsection{Preprocessing} \label{subsubsec:KnotPrePro}

First, the data was clipped to $132 \times 132$ locations. The baseline of the pressure-time course at each location was estimated by the median of the pre-excitation-pulse time points ($1$-$4$) and subtracted. Then, a zero-phase band-pass filter around $0.5$-$20$\si{\mega\hertz} was applied (see \verb|applyFilter.m| in the k-Wave toolbox). The $1$\% locations with the highest variance (computed by the time points $7$-$13$) were excluded from the further analysis. Finally, the remaining data was clipped to the time points $10$-$400$.

\subsubsection{Results} \label{subsubsec:KnotRes} 

\begin{figure}[tb]
  \centering
  \scriptsize
\begin{tabular*}{1\textwidth}{@{\extracolsep{-10pt}}cccc}
  \includegraphics[width=0.25\textwidth]{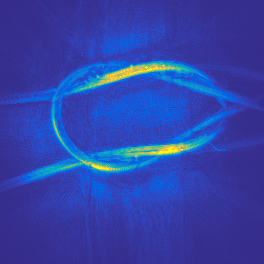} & 
 \includegraphics[width=0.25\textwidth]{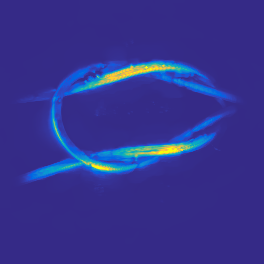} &
\includegraphics[width=0.25\textwidth]{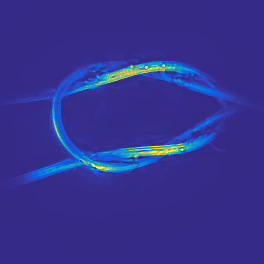} &
\includegraphics[width=0.25\textwidth]{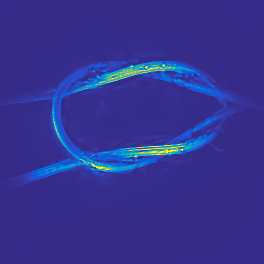} 
\\
  \includegraphics[width=0.25\textwidth]{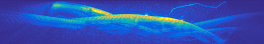} & 
 \includegraphics[width=0.25\textwidth]{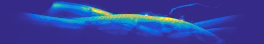} &
\includegraphics[width=0.25\textwidth]{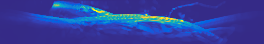} &
\includegraphics[width=0.25\textwidth]{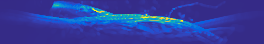} 
\\
  \includegraphics[width=0.25\textwidth]{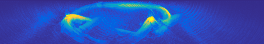} & 
 \includegraphics[width=0.25\textwidth]{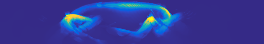} &
\includegraphics[width=0.25\textwidth]{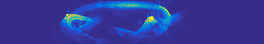} &
\includegraphics[width=0.25\textwidth]{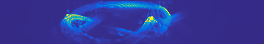} 
\\
(a) TR, cnv. & (b) TRppTV+, cnv. & (c) TV+, cnv. & (d) TV+Br, cnv. \\[4pt]
  \includegraphics[width=0.25\textwidth]{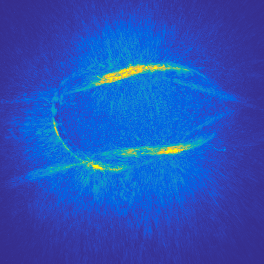} &
  \includegraphics[width=0.25\textwidth]{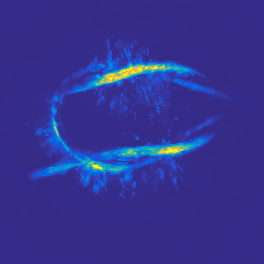} &
  \includegraphics[width=0.25\textwidth]{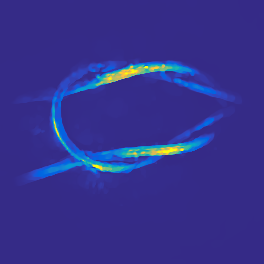} &
   \includegraphics[width=0.25\textwidth]{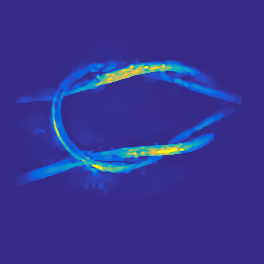} 
\\
  \includegraphics[width=0.25\textwidth]{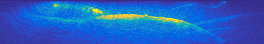} &
  \includegraphics[width=0.25\textwidth]{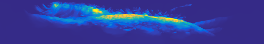} &
  \includegraphics[width=0.25\textwidth]{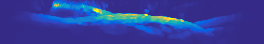} &
   \includegraphics[width=0.25\textwidth]{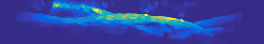} 
\\
  \includegraphics[width=0.25\textwidth]{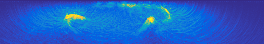} &
  \includegraphics[width=0.25\textwidth]{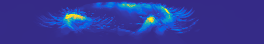} &
  \includegraphics[width=0.25\textwidth]{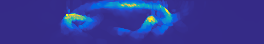} &
   \includegraphics[width=0.25\textwidth]{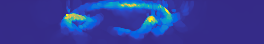} 
\\
(e) TR, rSP-4 & (f) TRppTV+, rSP-4 & (g) TV+, rSP-4 & (h) TV+Br, rSP-4\\[4pt]
  \includegraphics[width=0.25\textwidth]{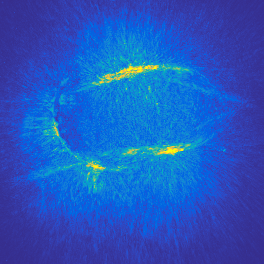} &
  \includegraphics[width=0.25\textwidth]{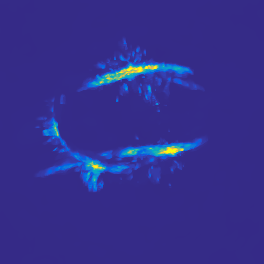} &
  \includegraphics[width=0.25\textwidth]{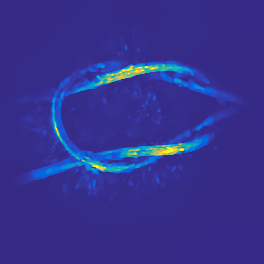} &
   \includegraphics[width=0.25\textwidth]{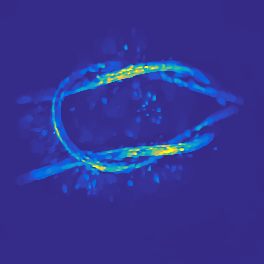} 
\\
  \includegraphics[width=0.25\textwidth]{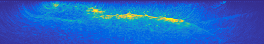} &
  \includegraphics[width=0.25\textwidth]{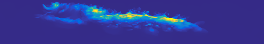} &
  \includegraphics[width=0.25\textwidth]{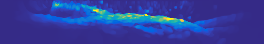} &
   \includegraphics[width=0.25\textwidth]{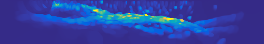} 
\\
  \includegraphics[width=0.25\textwidth]{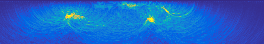} &
  \includegraphics[width=0.25\textwidth]{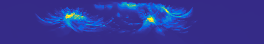} &
  \includegraphics[width=0.25\textwidth]{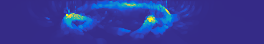} &
   \includegraphics[width=0.25\textwidth]{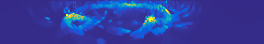} 
\\
(i) TR, rSP-8 & (j) TRppTV+, rSP-8 & (k) TV+, rSP-8 & (l) TV+Br, rSP-8\\[4pt]
  \includegraphics[width=0.25\textwidth]{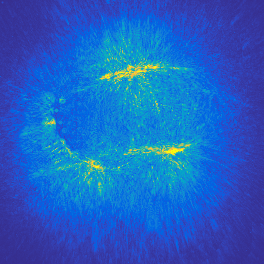} &
  \includegraphics[width=0.25\textwidth]{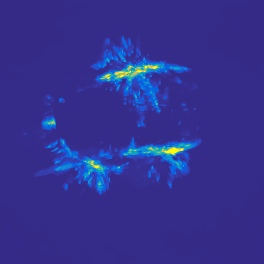} &
  \includegraphics[width=0.25\textwidth]{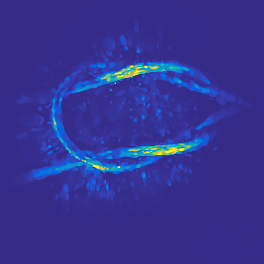} &
   \includegraphics[width=0.25\textwidth]{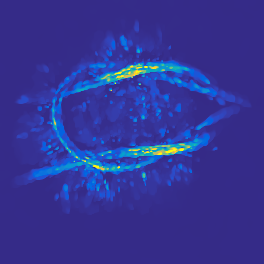} 
\\
  \includegraphics[width=0.25\textwidth]{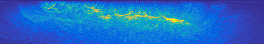} &
  \includegraphics[width=0.25\textwidth]{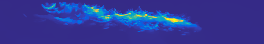} &
  \includegraphics[width=0.25\textwidth]{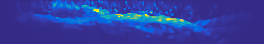} &
   \includegraphics[width=0.25\textwidth]{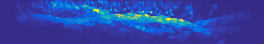} 
\\
  \includegraphics[width=0.25\textwidth]{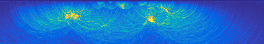} &
  \includegraphics[width=0.25\textwidth]{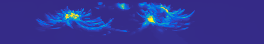} &
  \includegraphics[width=0.25\textwidth]{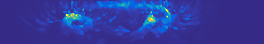} &
   \includegraphics[width=0.25\textwidth]{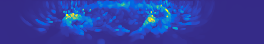} 
\\
(m) TR, rSP-16 & (n) TRppTV+, rSP-16 & (o) TV+, rSP-16 & (p) TV+Br, rSP-16
\end{tabular*}
\caption{Results for frame $23$ of the data set Knot: In each sub-figure, maximum intensity projections in X (top), Y (middle) and Z (bottom) direction are shown.}
\label{fig:ThKn}
\end{figure}

Table \ref{tbl:ModelPara} shows the parameters of the acoustic model used for the inversion. Note that the spacing of the spatial grid, $\Delta_{x/y/z}$, is $2$ times finer than the distance between scanning locations: Assuming a sound speed of $1540$ \si{\meter \per \second} the Nyquist criterion would require that we had sampled with $\delta_{y/z} = 38.5$ \si{\micro \meter} in space and $\delta_t = 25$ \si{\nano\second} in time to be able to reconstruct initial pressure distributions leading to signals with a frequency content of $20$ \si{\mega \hertz}. This means that the conventional data is over-sampled in time, but already under-sampled in space, similar to the simulation studies with Tumor2 (cf. Section  \ref{subsec:BackNyquist}). When choosing a finer spatial grid spacing to reconstruct the data as compared to the one in which it was recorded we attempt to recover some spatial resolution from the higher temporal resolution of the pressure time series. \\
In 4D PAT, we can vary the sub-sampling operator $C$ used in each frame $i$. For single-point sub-sampling, one would try to avoid measuring the pressure time series at the same location in subsequent frames as it may contain very similar information\footnote{For patterned interrogation, one would not use the same pattern $\phi_j$ in subsequent frames}. Therefore, we randomly partitioning the set of all scanned locations into $M_{sub}$ subsets, each containing $M_c = M/M_{sub}$ different random locations. This yields a sequence of $M_{sub}$ sub-sampling operators $C_i$, $i=1,\ldots,M_{sub}$ that we periodically apply to the set of all $45$ frames, i.e., after $M_{sub}$ frames, each locations has been scanned once and the $M_{sub}+1$-th frame is scanned with $C_1$, again. Figure \ref{fig:ThKn} shows the inverse reconstructions of the middle frame $23$ for conventional data, rSP-4, rSP-8 and rSP-16 (a movie of the complete reconstruction can be found in the supplementary material). Next to TR, TV+ and TV+Br, we also show the result of post-processing the TR solution $p_{\rm TR}$ with positivity-constrained TV denoising, which we will denote by "TRppTV+":
\begin{equation}
p_{\text{\tiny TRppTV+}} \mydef \argminsub{p \geqslant 0} \left\lbrace \frac{1}{2} \sqnorm{p - p_{\rm TR}} +  \lambda_{\rm pp} \text{TV}(p) \right\rbrace, \label{eq:TV+Denoise}
\end{equation}
where we chose $\lambda_{\rm pp}$ large enough to suppress most visible reconstruction artifacts. Solving \eref{eq:TV+Denoise} is discussed in \ref{sec:Opti}. The regularization parameter chosen for TV+ was $\lambda_{\text{\tiny TV+}}$ = $0.02$ for the conventional data while this value was multiplied by $4/M_{sub}$ for the sub-sampled data. For TV+Br, we carried out 10 Bregman iterations with $\lambda_{\text{\tiny TV+Br}}$ = $12.5 \lambda_{\text{\tiny TV+}}$.

\subsection{Patterned Interrogation - Static Phantom} \label{subsec:HrKn}

Next, we investigate data acquired by the patterned interrogation FP scanner (cf. Figure \ref{subfig:FPBP}), measuring a static experimental phantom, which we will call "Hair".

\subsubsection{Setup} \label{subsubsec:HrKnSetup}

The full technical details of the scan can be found in \cite{HuZhBeArBeCo14}. The target to be scanned was a knotted artificial hair (see Figure \ref{subfig:HrKnSetup}, diameter $\sim$150\si{\micro\meter} ), immersed in 1\% Intralipid solution and positioned approximately $2$\si{\milli\meter} above the detection plane and $3$\si{\milli\meter} under the Intralipid surface. On the DMD, an active area of $640 \times 640$ micromirrors was subdivided by grouping $5 \times 5$ micromirrors to form one of $128 \times 128$ "pixels". Due to an angle in the optical path, each pixel corresponded to an area of $62.12$ $\times$ $68$ \si{\micro\meter} on the detection plane. Then, the rows of a $16384 \times 16384$ scrambled Hadamard matrix were used to implement $128 \times 128$ binary pattern on the DMD.

\subsubsection{Preprocessing} \label{subsubsec:HrKnPrePro}

First, the data needed to be calibrated to fit to our model (cf. Section \ref{subsec:DiscFwdModel}). Subsequently, a zero-phase band-pass filter around $3$-$22.5$\si{\mega\hertz} was applied and the data was clipped to the time points $85$-$140$.

\subsubsection{Results} \label{subsubsec:HrKnRes}

\begin{figure}[tb]
  \centering
  \scriptsize
\begin{tabular*}{1\textwidth}{@{\extracolsep{-10pt}}ccccc}
  \includegraphics[width=0.2\textwidth]{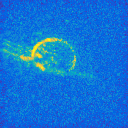} & \includegraphics[width=0.2\textwidth]{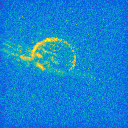} & 
  \includegraphics[width=0.2\textwidth]{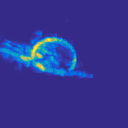} & 
  \includegraphics[width=0.2\textwidth]{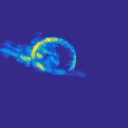} & 
  \includegraphics[width=0.2\textwidth]{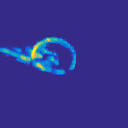} \\
  \includegraphics[width=0.2\textwidth]{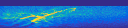} &  \includegraphics[width=0.2\textwidth]{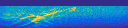} &  \includegraphics[width=0.2\textwidth]{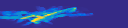} &  \includegraphics[width=0.2\textwidth]{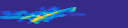} & \includegraphics[width=0.2\textwidth]{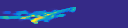} \\
  \includegraphics[width=0.2\textwidth]{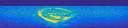} & \includegraphics[width=0.2\textwidth]{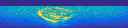}  &  \includegraphics[width=0.2\textwidth]{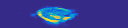} &  \includegraphics[width=0.2\textwidth]{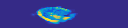} & \includegraphics[width=0.2\textwidth]{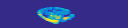} \\
(a) TR, sHd-1 & (b) BP, sHd-1 & (c) TRppTV+, sHd-1 & (d) BPppTV+, sHd-1 & (e) TV+Br, sHd-1\\[6pt]
  \includegraphics[width=0.2\textwidth]{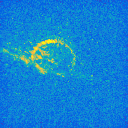} & \includegraphics[width=0.2\textwidth]{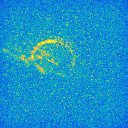}  & \includegraphics[width=0.2\textwidth]{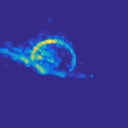} & \includegraphics[width=0.2\textwidth]{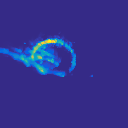}  & \includegraphics[width=0.2\textwidth]{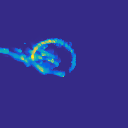} \\
  \includegraphics[width=0.2\textwidth]{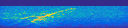} & \includegraphics[width=0.2\textwidth]{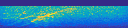}  & \includegraphics[width=0.2\textwidth]{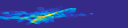} & \includegraphics[width=0.2\textwidth]{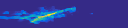} & \includegraphics[width=0.2\textwidth]{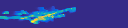} \\
  \includegraphics[width=0.2\textwidth]{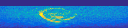} & \includegraphics[width=0.2\textwidth]{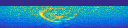}  & \includegraphics[width=0.2\textwidth]{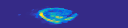} & \includegraphics[width=0.2\textwidth]{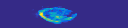} &\includegraphics[width=0.2\textwidth]{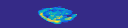} \\
(f) TR, sHd-4 & (g) BP, sHd-4 & (h) TRppTV+, sHd-4 & (i) BPppTV+, sHd-4 & (j) TV+Br, sHd-4\\
  \includegraphics[width=0.2\textwidth]{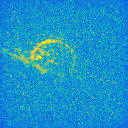} & \includegraphics[width=0.2\textwidth]{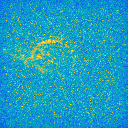}  &  \includegraphics[width=0.2\textwidth]{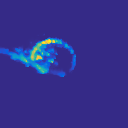} & \includegraphics[width=0.2\textwidth]{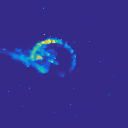} & \includegraphics[width=0.2\textwidth]{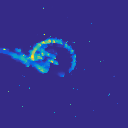} \\
  \includegraphics[width=0.2\textwidth]{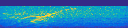} & \includegraphics[width=0.2\textwidth]{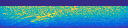} &  \includegraphics[width=0.2\textwidth]{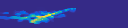} & \includegraphics[width=0.2\textwidth]{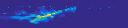} & \includegraphics[width=0.2\textwidth]{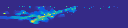} \\
  \includegraphics[width=0.2\textwidth]{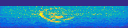} & \includegraphics[width=0.2\textwidth]{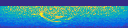} &  \includegraphics[width=0.2\textwidth]{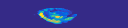} & \includegraphics[width=0.2\textwidth]{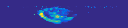} & \includegraphics[width=0.2\textwidth]{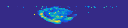} \\
(k) TR, sHd-8 & (l) BP, sHd-8 & (m) TRppTV+, sHd-8 & (n) BppTV+, sHd-8 & (o) TV+Br sHd-8\\
  \includegraphics[width=0.2\textwidth]{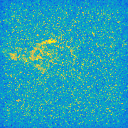} & \includegraphics[width=0.2\textwidth]{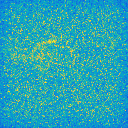} &  
  \includegraphics[width=0.2\textwidth]{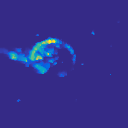} & \includegraphics[width=0.2\textwidth]{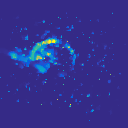} & \includegraphics[width=0.2\textwidth]{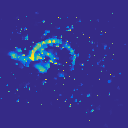} \\
  \includegraphics[width=0.2\textwidth]{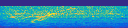} & \includegraphics[width=0.2\textwidth]{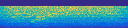} &  
  \includegraphics[width=0.2\textwidth]{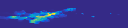} & \includegraphics[width=0.2\textwidth]{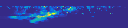} & \includegraphics[width=0.2\textwidth]{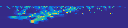} \\
  \includegraphics[width=0.2\textwidth]{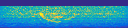} & \includegraphics[width=0.2\textwidth]{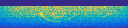} &  
  \includegraphics[width=0.2\textwidth]{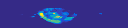} & \includegraphics[width=0.2\textwidth]{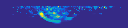} & \includegraphics[width=0.2\textwidth]{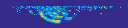} \\
(p) TR, sHd-16 & (q) BP, sHd-16 & (r) TRppTV+, sHd-16 & (s) BPppTV+, sHd-16  & (t) TV+Br, sHd-16
\end{tabular*}
\caption{Results for data set Hair: In each sub-figure, maximum intensity projections in X (top), Y (middle) and Z (bottom) direction are shown.}
\label{fig:HrKn}
\end{figure}

Table \ref{tbl:ModelPara} shows the parameters of the acoustic model used for the inversion. Note that the conventional data is, again, slightly over-sampled in time but under-sampled in space. For all computed solutions $p$, all voxels in the first $6$ x-layers were set to $0$ in a post-processing step. The latter was done to ease the visualization of the results through maximum intensity projections: The first $12$ time points of the signal only seem to contain noise, which means that up to a distance of $12 c_0 \Delta_t/\Delta_x = 5.8$ in voxel length, $p_0$ will only account for noise. For TR and BP solutions, voxels with negative values were also set to $0$. Figure \ref{fig:HrKn} shows different reconstructions using the conventional data (sHd-1) and $M_{sub} = 4,8,16$: The regularization parameter chosen for TV+Br was $\lambda_{\text{\tiny TV+Br}}$ = $1.5 \cdot 10^{-4}$ for the conventional data while this value was first divided by $M_{sub}$ and then multiplied by $1/1.2/1.4$ for $M_{sub} = 4/8/16$ for the sub-sampled data. A total of 10 Bregman iterations were carried out.

\subsection{In-Vivo Measurements - Single Point Sub-Sampling} \label{subsec:BlVe}

While validating inverse methods on data from experimental phantoms is an important step forward from pure simulation studies, experimental phantoms cannot reproduce all the features of real in-vivo data sets. As a last example, we therefore investigate a static in-vivo data set of skin vasculature and subcutaneous anatomy near the right flank of a nude mouse, which we will call "Vessels".

\subsubsection{Setup} \label{subsubsec:BlVeSetup}

The data was acquired with an excitation wavelength of $590$\si{\nano\meter}, further technical details and illustrations can be found in  \cite{JaLaOgTrCoZhJoPiPhMaLyPePuBe15}. A conventional scan consisted of $142 \times 141$ locations over a region of size of $14$\si{\milli\meter} $\times$ $14$\si{\milli\meter} ($\delta_x = \delta_y = 100$\si{\micro\meter}), measured at $M_t = 630$ time points with a resolution of $\delta_t = 10$\si{\nano\second}.

\subsubsection{Preprocessing} \label{subsubsec:BlVePrePro}

The data was clipped to $141 \times 141$ locations and to the time points $10$-$630$. 
         
\subsubsection{Results} \label{subsubsec:BlVeRes}

\begin{figure}[tb]
  \centering
  \scriptsize
\begin{tabular*}{1\textwidth}{@{\extracolsep{-10pt}}ccc}
\includegraphics[width=0.33\textwidth]{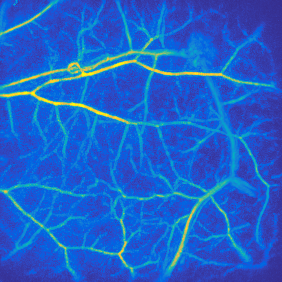} &   
  \includegraphics[width=0.33\textwidth]{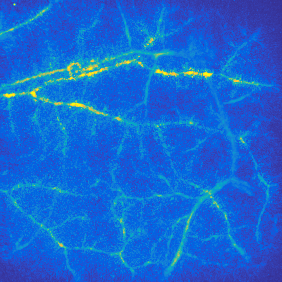} &   \includegraphics[width=0.33\textwidth]{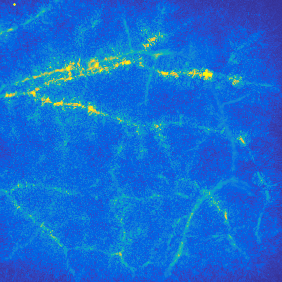} \\
  \includegraphics[width=0.33\textwidth]{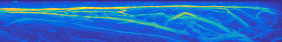} &   
  \includegraphics[width=0.33\textwidth]{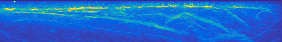} &   \includegraphics[width=0.33\textwidth]{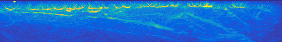} \\
  \includegraphics[width=0.33\textwidth]{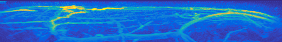} &   
  \includegraphics[width=0.33\textwidth]{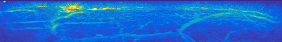} &   \includegraphics[width=0.33\textwidth]{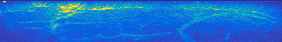} \\
\includegraphics[width=0.33\textwidth]{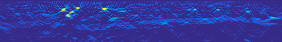} & \includegraphics[width=0.33\textwidth]{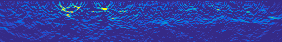} & \includegraphics[width=0.33\textwidth]{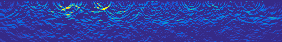} \\
(a) TR, cnv. & (b) TR, rSP-4 & (c) TR, rSP-8  \\[6pt]
  \includegraphics[width=0.33\textwidth]{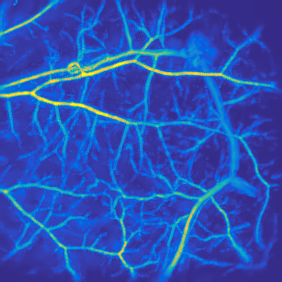} &   
  \includegraphics[width=0.33\textwidth]{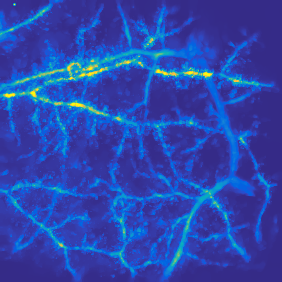} &   \includegraphics[width=0.33\textwidth]{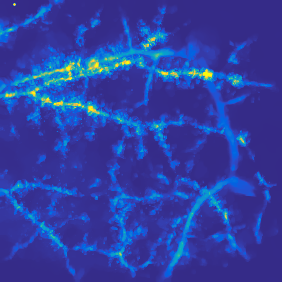} \\
  \includegraphics[width=0.33\textwidth]{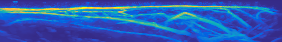} &   
  \includegraphics[width=0.33\textwidth]{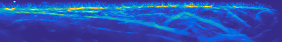} &   \includegraphics[width=0.33\textwidth]{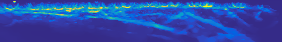} \\
  \includegraphics[width=0.33\textwidth]{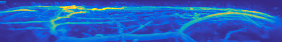} &   
  \includegraphics[width=0.33\textwidth]{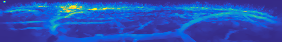} &   \includegraphics[width=0.33\textwidth]{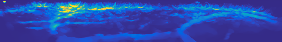} \\
\includegraphics[width=0.33\textwidth]{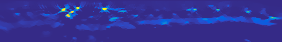} & \includegraphics[width=0.33\textwidth]{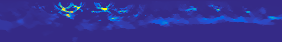} & \includegraphics[width=0.33\textwidth]{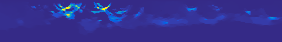} \\
(d) TRppTv+, cnv. & (e) TRppTV+, rSP-4 & (f) TRppTv+, rSP-8  \\[6pt]
\includegraphics[width=0.33\textwidth]{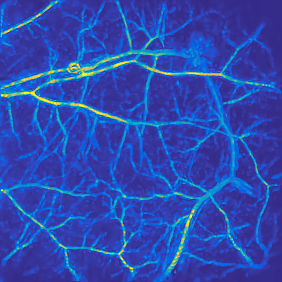} & \includegraphics[width=0.33\textwidth]{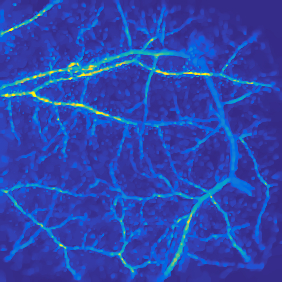} & \includegraphics[width=0.33\textwidth]{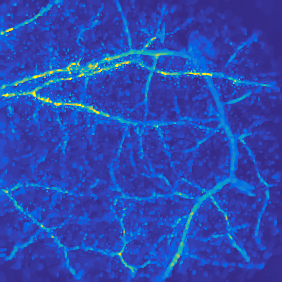}  \\
\includegraphics[width=0.33\textwidth]{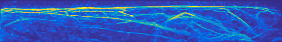} & \includegraphics[width=0.33\textwidth]{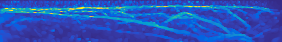} & \includegraphics[width=0.33\textwidth]{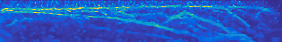}  \\
\includegraphics[width=0.33\textwidth]{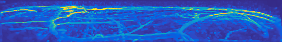} & \includegraphics[width=0.33\textwidth]{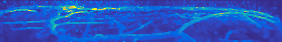} & \includegraphics[width=0.33\textwidth]{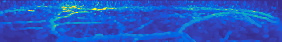}  \\
\includegraphics[width=0.33\textwidth]{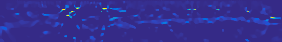} &
\includegraphics[width=0.33\textwidth]{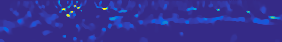} & \includegraphics[width=0.33\textwidth]{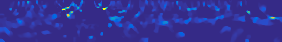} \\
(g) TV+Br, cnv. & (h) TV+Br, rSP-4 & (i) TV+Br, rSP-8  
\end{tabular*}
\caption{Results for in-vivo data set Vessels: In each sub-figure, from top to bottom:  Maximum intensity projections in X, Y, Z direction and slice through $z=74$.}
\label{fig:BlVe}
\end{figure}

Table \ref{tbl:ModelPara} shows the parameters of the acoustic model used for the inversion. Note that by a similar reasoning about the spatial and temporal sampling intervals, the spatial spacing $\Delta_{x/y/z}$ is, again, chosen $2$ times finer than the distance between scanning locations. Figure \ref{fig:BlVe} shows maximum intensity projections and a slice through $z=74$ for TR, TRppTV+ and TV+Br solutions when using the conventional, rSP-4 and rSP-8 data. The denoising parameter $\lambda_{\rm pp}$ , was, again, chosen large enough to suppress most visible reconstruction artifacts. The regularization parameter chosen for TV+Br was $\lambda_{\text{\tiny TV+Br}}$ = $6.25 \cdot 10^{-2}$ for the conventional data while this value was multiplied by $4/M_{sub}$ for the sub-sampled data. A total of 10 Bregman iterations were carried out. 

\section{Discussion, Outlook and Conclusion}\label{sec:DisOutCon}

\subsection{Discussion} \label{subsec:Dis}

The main results of the simulation studies (Section \ref{sec:SimStud}) can be summarized as follows:
\begin{itemize}
\item Using model-based variational reconstruction methods employing spatial sparsity constraints, such as TV+, \eref{eq:TV+}, is essential for obtaining good quality PA images from sub-sampled PAT data. The linear methods, TR \eref{eq:TRnum} and BP \eref{eq:BP}, and the L2+ method \eref{eq:L2+} could not produce images of acceptable quality in any setting.
\item The results obtained for Tumor1 and Tumor2 demonstrate that the image quality obtained for a certain sub-sampling rate $M_{sub}$ varies strongly: The "inverse crime" data of the more superficial, high contrast target Tumor1 can be up-sampled up to $M_{sub} = 128$ without a significant loss of image quality. On the other hand, a bad \emph{model-fit}, i.e., a mismatch between of the models used for data generation and inversion, combined with a more challenging target, such as Tumor2, significantly impairs the image quality beyond a certain sub-sampling rate ($M_{sub} = 16-32$ in our particular example). 
\item Using Bregman iterations improves upon conventional variational approaches for PAT. Most importantly, the systematic contrast loss of small scale structures such as blood vessels is mitigated which is a crucial prerequisite for QPAT studies.
\item Sub-sampling by patterned interrogation (sHd) was slightly more efficient than single point sub-sampling (rSP). 
\end{itemize}
The sub-sampling rates we achieved with experimental data were similar to those obtained with simulated data in the more realistic Tumor2 scenario. However, unexpectedly, the qualitative difference between the linear methods TR and BP and the variational methods TV+ and TV+Br was not as dramatic as in the simulation studies (cf. Figures \ref{fig:EasyLinRes}, \ref{fig:EasyVarRes}, \ref{fig:ThKn}, \ref{fig:HrKn}). In many cases, post-processing the TR solution with TV+ denoising, \eref{eq:TV+Denoise}, comes remarkably close to the TV+ solution, which is not to be expected from the simulation studies (see, e.g., the TR sHd-128 solution in Figure \ref{fig:EasyLinRes}). In addition, Bregman iterations could not improve much over conventional variational approaches, again in contrast to our findings in Section \ref{subsubsec:EnhanceBregIter}. A partial explanation for both findings is a bad model fit: While some pre-processing routines (e.g., baseline correction, band-pass filtering and a detection and deletion of corrupt channels) were implemented and carried out to align the data with the model used, other known model-mismatches (e.g., the inhomogeneous sensitivity of the FP sensor and the non-whiteness of the measurement noise) were not accounted for. \\
A closer examination of the "Knot" data reveals several flaws which deteriorated our results: Firstly, the optical excitation was inhomogeneous in lateral direction (cf. Figure \ref{fig:ThKn}), leading to an inhomogeneous initial pressure distribution in regions consisting of the same materials. The TV energy is not well-suited to recover such targets. Secondly, a the baseline shifts in the data are more complex than what we corrected for and a spatio-temporal visualization of the data shows several artifacts on the sensor that our automatic channel deletion procedure cannot fully remove. And lastly, the acoustic properties of the polythene tubes lead to reflections we did not account for in our model. For these reasons, the decrease in image quality for sub-sampled data is faster compared to the simulation studies. While we see a clear advantage of using TV+ as opposed to the simpler TRppTV+, using Bregman iterations does not seem to lead to a better image quality. \\
The "Hair" data was acquired with the novel patterned FP scanner which still suffers from several major technical difficulties, e.g., the widening of the interrogation beam leads to a significant loss in SNR and systematic shifts in signal power can be observed that need to be examined more carefully. Furthermore, we suspect that the hair knot might have moved during the acquisition (cf Figure \ref{fig:HrKn}(c)). For these reasons, already the reconstructions from the sHd-1 data are not of very good quality. If we accept these as a ground truth nonetheless, we see that we can reach sub-sampling rates of around $M_{sub} = 8$ without a significant further loss of image quality. If we compare the different methods, we find that TRppTV+ yields the visually most appealing results while BPppTV+ and TV+Br look very similar. As BPppTV+ is more or less the first iterate of the optimization scheme we use to compute TV+Br (cf. Section \ref{sec:Opti}), the latter might, again, point to a bad model-fit (see above). \\
For the "Vessel" data set, the visual impression of the conventional data reconstructions (cf. Figure \ref{fig:BlVe}) and the fact that we did not have to perform any pre-processing suggest that the in-vivo data examined is of good quality and we have a good model-fit. A comparison with the results from experimental phantoms further reveals that the diffusive nature of biological tissue leads to a more even lateral illumination. Now, TV+Br clearly outperforms TR and TRppTV+. For instance, from the slice view we can see that TR seems to overestimate the diameter of the blood vessels. The reason for not achieving high sub-sampling rates despite the good data quality is the apparent mismatch between the target geometry and the spatial sparsity constraints employed by the TV energy: The PA image is exceptionally rich in vasculature, but TV regularization tends to break up such anisotropic, line-like structures (cf. Figure \ref{fig:BlVe}(i)).

\subsection{Outlook}\label{subsec:Out}

As the simulation studies show that a model misfit can severely decrease the sub-sampling rates achievable, we need to improve the accuracy of the acoustic forward model to obtain better results for experimental data: Data pre-processing aligns the data with the forward model, \emph{model calibration} can determine some uncertain parameters of the forward model (such as the FP sensitivity distribution or the noise statistics), and Bayesian \emph{model selection} \cite{To11} or Bayesian \emph{approximation error modeling} \cite{ArKaKoScSoTaVa06,KaSo07,TaPuCoKaAr13} can reduce or account for the uncertainty in other parameters, such as $c_0$. To improve the results for in-vivo data (cf. Figure \ref{fig:BlVe}) acoustic absorption models of biological tissue \cite{TrZhCo10,Tr13} need to be incorporated.\\
While we exploited spatial sparsity to accelerate the acquisition of a single scan here, the next step to enable 4D PAT imaging with both high spatial and temporal resolution (cf. Section  \ref{subsec:Knot}) is to extend the frame-by-frame inversion methods examined here to \emph{full spatio-temporal} variational models that also exploit the temporal redundancy of data generated by dynamics of low complexity.\\
We used the TV energy as a generic, well-understood first example of a spatial sparsity constraint. However, as discussed above, it is, e.g., not very suitable to recover thin vasculature. Higher sub-sampling rates could be reached by employing more sophisticated regularization functionals designed to recover such anisotropic structures \cite{KuLi12}.\\
Choosing random locations as single-point sub-sampling and scrambled Hadamard pattern as patterned interrogation was based on results obtained for similar applications such as CT and MRI. The optimal choice for PAT applications is yet to be determined. For instance, \cite{ScGlZaLiDrSc15} has shown that a non-uniform distribution of the sampling locations in single-point sub-sampling can be used to focus into a specific area (at the expense of the resolution elsewhere). A theoretical examination, e.g., through \emph{micro-local analysis} \cite{FrQu15}, could help to gain new insights on this.


\subsection{Conclusion}\label{subsec:Con}

In this study, we investigated different possibilities to sub-sample the incident photoacoustic field in order to accelerate the acquisition of high resolution PAT. In simulation studies, we demonstrated that PAT wave fields generated by targets with a low spatial complexity can indeed be highly compressible and identified under which conditions this feature can be exploited to obtain high quality images from highly sub-sampled data: Firstly, variational image reconstruction methods employing sparsity constraints that match the structure of the target have to be used. Secondly, using an accurate forward model well-aligned with the data is crucial. We furthermore applied the methods developed to three experimental data sets from experimental phantoms and in-vivo recordings. While obtaining promising first results, we also identified several challenges for realizing the full potential of data sub-sampling, most notably obtaining a good model-fit as discussed above. While we focused on sub-sampling the data using the Fabry-P\'{e}rot based scanner, the novel reconstruction strategies offer new opportunities to dramatically increase the acquisition speed of other photoacoustic scanners that employ point-by-point sequential scanning as well as reducing the channel count of parallelized schemes that use detector arrays.


\ack

We gratefully acknowledge the support of NVIDIA Corporation with the donation of the Tesla K40 GPU used for this research.\\
The numerical phantom is based upon \si{\micro}CT data that was kindly provided by Simon Walker-Samuel at the Centre for Advanced Biomedical Imaging, University College London, and a segmentation thereof kindly provided by Roman Hochuli, University College London.\\
This work was supported by the Engineering and Physical Sciences Research Council, UK (EP/K009745/1), the European Union project FAMOS (FP7 ICT, Contract 317744) and the National Institute of General Medical Sciences of the National Institutes of Health under grant number P41 GM103545-18.

\appendix

\section{Discrete Total Variation Energy} \label{sec:TV}

Let the voxels of the 3D pressure $p \in \R^N$, with $N = N_x N_y N_z$ be indexed by $(i,j,k)$, $i = 1,\ldots,N_x$, $j = 1,\ldots,N_y$, $k=1,\ldots,N_z$. Using finite forward differences, the most commonly used discretization of the total variation seminorm with Neumann boundary conditions is given by   
\begin{equation*}
\fl \text{TV}(p) = \sum_{(i,j,k)} \sqrt{(p_{(i+1,j,k)}-p_{(i,j,k)})^2 + (p_{(i,j+1,k)} - p_{(i,j,k)})^2 + (p_{(i,j,k+1)} - p_{(i,j,k)})^2}, 
\end{equation*}
where $p_{(N_x+1,j,k)} := p_{(N_x,j,k)}$, $p_{(i,N_y+1,k)} := p_{(i,N_y,k)}$ and $p_{(i,j,N_z+1)} := p_{(i,j,N_z)}$. Additional terms have to be added to account for different boundary conditions: In the simulation studies, we use Dirichlet boundary conditions which requires to add the jumps at the domain boundaries. For the experimental data, we impose Dirichlet boundary conditions at the detection plane voxels $x = 0$ while Neumann boundary conditions are applied on all other faces of the image cube.


\section{Optimization}\label{sec:Opti}

The optimization problem given by \eref{eq:VarReg} consists of the minimizing the sum of two functionals $\mathcal{E}(p) =  \mathcal{D}(p) + \lambda \mathcal{J}(p)$, where we can compute the gradient of the strictly convex, smooth functional $\mathcal{D}(p)$,
\begin{equation}
\mathcal{D}(p) = \frac{1}{2} \sqnorm{C A \, p - f^c}, \qquad \nabla \mathcal{D}(p) = A^T C^T \left(C A \, p - f^c \right),
\end{equation}
but no higher derivatives, and we know how to compute the \emph{proximal operator} of the convex, potentially non-differentiable functional $\mathcal{J}(p)$:
\begin{equation}
\text{prox}_{\mathcal{J},\alpha}(p) \mydef \argminsub{q} \left\lbrace \J(q) + \frac{1}{2 \alpha} \sqnorm{q - p} \right\rbrace \label{eq:Prox}
\end{equation}
In the case of $\J(p)$ being the positivity-constrained TV energy, the proximal operators simply solves a positivity-constrained TV denoising problem \eref{eq:TV+Denoise}. \\
A wide range of \emph{semi-smooth, first order optimization} algorithms for image reconstruction have been developed over recent years \cite{BuSaSt14}, each of them advantageous for a specific scenario. In our case, the special feature of PAT is that the application of $A$ and $A^T$ requires considerably more computation time than solving most proximal operators \eref{eq:Prox}, including the 3D positivity-constrained TV denoising problem \eref{eq:TV+Denoise}, up to a high numerical precision. Under these circumstances, the rather simple \emph{proximal gradient descent} scheme:
\begin{equation}
\fl \quad p^{k+1} = \text{prox}_{\J,\eta \lambda} \left(p^k - \eta A^T C^T \left(C A p^k -f^c \right) \right), \qquad p^0 = 0, \; k=1,\ldots,K \label{eq:ProxGradDesc}
\end{equation}
turns out to be most efficient if tuned carefully (see \cite{GoStBa14} for an extensive overview):
\begin{itemize}
  \item The step-size $\eta$ is set to $1.8/L$, where $L$ is an approximation of the Lipschitz constant of $A^T C^T C A$. For a given setting and sub-sampling scheme, $L$ can be computed quite efficiently by a simple power iteration and then stored in a look-up table.  
  \item We use a gradient extrapolation modification (\emph{fast or accelerated gradient methods}) ensuring a quadratic convergence. The concrete technique we use is the \emph{FISTA} algorithm \cite{BeTe09}, where we restart the acceleration if an increase in $\mathcal{E}(p^k)$ is detected and switch to a normal gradient for this iterations $k$, followed by up to 5 backtracking steps if necessary ($\eta$ is, however, not changed for future iterations).  
  \item For the 3D PAT problems considered in this work, the iterates from $k = 11$ onwards are usually visually not distinguishable. However, we computed a maximum of $K = 50$ iterations for all except for Knot, where we only computed $K = 20$ iterations. We terminated the iteration earlier if $p^k$ did not change or $\min_k \mathcal{E}(p^k)$ did not decrease for $5$ times in a row.
\end{itemize}
The proximal operator for \eref{eq:L2+} can be computed component-wise and explicitly:
\begin{equation}
\fl \qquad \tilde{p} = \argminsub{q \geqslant 0} \left\lbrace \sqnorm{q} + \frac{1}{2 \alpha} \sqnorm{q - p} \right\rbrace \quad \Leftrightarrow \quad  \tilde{p}_i = \max \left(0,\frac{q_i}{1+\alpha} \right)
\end{equation}
The positivity-constrained TV denoising is implemented by a \emph{primal-dual hybrid gradient} algorithm as described in \cite{ChPo11}.


\section{Implementation}\label{sec:Implementation}

All routines have been implemented in as part of a larger, Matlab toolbox for PAT image reconstruction which will be made available in near future. The toolbox relies on the k-Wave toolbox (see \cite{TrCo10},  \href{http://www.k-wave.org/}{http://www.k-wave.org/}) to implement $A$ and $A^T$, which allows to use highly optimized C++ and CUDA code to compute the 3D wave propagation on parallel CPU or GPU architectures. To give an idea about the range of different computations times, computing one application of $A$ for the in-vivo Vessels scenario (cf. Table \ref{tbl:ModelPara}) in single precision takes 15\si{\second} using the optimized CUDA code on a Tesla K40 GPU (counting only the GPU run-time), 51\si{\second} using the Matlab code on the same GPU, 47\si{\second}/6\si{\minute} 36\si{\second} using the optimized C++ code on 12/1 cores of an Intel Xeon CPU (2.70GHz) (counting only the CPU run-time) and 4\si{\minute} 3\si{\second}/26\si{\minute} 48\si{\second} using the Matlab code on 12/1 cores of the same CPU.


\section*{References}

\bibliographystyle{unsrt}
\bibliography{all}


\noappendix

\Tables

\begin{table}
\caption{\label{tbl:Abb} List of commonly occurring abbreviations.} 
\lineup
\begin{tabular}{@{}lll}        
\br                     
Abbreviation & Meaning & Reference  \\
\mr
BP  & \textit{back projection} & Sec. \ref{subsec:TwoStep}, \eref{eq:BP} \\
BPppTV+  & BP followed by TV+ denoising as post-processing & see TRppTV+ \\
cnv.  &  conventional (data sampling) & see Sec. \ref{subsubsec:SiPo} \\
DMD  &  \textit{digital micromirror device} &  Sec. \ref{subsec:DiscFwdModel} \\
DP  &  \textit{discrepancy principle} & Sec. \ref{subsec:SimVaTu1}, \eref{eq:DiscPr} \\
FP  & \textit{Fabry-P\'{e}rot} interferometer & Sec. \ref{subsec:CSFP} \\
L2+  & positivity-constrained $\ell_2$ regularization & Sec. \ref{subsec:VarImRec}, \eref{eq:L2+} \\
mxIP & \textit{maximum intensity projection} & Sec. \ref{subsec:NumPhan} \\
PSNR & \textit{peak signal-to-noise ratio} & \eref{eq:PSNR} \\
(Q)PAT & \textit{(quantitative)} photoacoustic tomography & Sec. \ref{sec:Intro} \\
rSP  & \textit{random single point} sub-sampling & Sec. \ref{subsubsec:SiPo}\\
sHd  & \textit{scrambled Hadamard} sub-sampling & Sec. \ref{subsubsec:PatInter}, \ref{subsec:DiscFwdModel} \\
TR  & \textit{time reversal} & Sec. \ref{subsec:TwoStep}, \eref{eq:TRnum} \\
TRppTV+  & TR followed by TV+ denoising as post-processing & \eref{eq:TV+Denoise} \\
TV+  & positivity-constrained total variation regularization & Sec. \ref{subsec:VarImRec},  \eref{eq:TV+} \\
TV+BR  & Bregman iterations applied to TV+ & Sec. \ref{subsec:BregIter},  \eref{eq:BregIterA}, \eref{eq:BregIterB} \\
\br
\end{tabular}
\end{table}

\begin{table}
\caption{\label{tbl:ModelPara} Parameter of the different inversion models used.} 
\lineup
\begin{tabular}{@{}lcccc}        
\br                     
parameter & Tumor1/2 & Knot & Hair & Vessels \\
\mr
$(N_x,N_y,N_z)$  &  128 & (44,264,264)   & (28,128,128)  & (42,282,282)     \\
$(\Delta_x,\Delta_y,\Delta_z)$ $[$\si{\micro\meter}$]$  & 156.25  &  75  & (62.12,62.12,68.00)  & 50   \\
$M_t$  &  740 & 391  & 56  & 621   \\
$\Delta_t/\delta_t$ $[$\si{\nano\second}$]$  & 31.25 &  12  & 20  &  10  \\
$M$ & 16384 & 17424 & 16384 & 19881 \\
$(\delta_x,\delta_y)$  & 156.25/312.5 &  150  & /  &  100  \\
$c_0$ $[$\si{\meter\per\second}$]$ & 1500  &  1540  & 1500 &  1420 \\
\br
\end{tabular}
\end{table}

\end{document}